\documentclass[a4paper,11pt]{article}
\usepackage{graphicx}
\usepackage{epsfig}
\usepackage{caption}
\usepackage{subcaption}
\usepackage{amsmath}
\usepackage{amssymb}
\usepackage{mathtools}
\usepackage[top=3cm, bottom=3cm, left=1.5cm, right=1.5cm]{geometry}
\usepackage{microtype}
\usepackage{lmodern}
\usepackage[utf8]{inputenc}
\usepackage{color}
\usepackage{bera}

\DeclarePairedDelimiter{\ceil}{\lceil}{\rceil}

\newcommand{\ds}{\displaystyle }

\newcommand{\beq}{\begin{equation} }
\newcommand{\eeq}{\end{equation}}
\author{{\Large Thomas M. Michelitsch$^a$, Federico Polito$^b$ and
Alejandro P. Riascos$^c$ }\\ \\
\footnotesize{$^a$ Sorbonne Universit\'e, Institut Jean le Rond d’Alembert,
CNRS UMR 7190} \\
\footnotesize{4 place Jussieu, 75252 Paris cedex 05, France} \\
\footnotesize{E-mail: michel@lmm.jussieu.fr}\\ 
\footnotesize{ORCID-ID: 0000-0001-7955-6666} \\ [1ex]
\footnotesize{$^b$ Department of Mathematics ``Giuseppe Peano'', University of Torino, Italy}  \\
\footnotesize{E-mail: federico.polito@unito.it} \\
\footnotesize{ORCID-ID: 0000-0003-1971-214X}\\ [1ex]
\footnotesize{$^c$Instituto de F\'isica, Universidad Nacional Aut\'onoma de M\'exico} \\
\footnotesize{Apartado Postal 20-364, 01000 Ciudad de M\'exico, M\'exico}  \\
\footnotesize{E-mail: aperezr@fisica.unam.mx} \\
\footnotesize{ORCID-ID: 0000-0002-9243-3246}
}
\title{Asymmetric random walks with bias generated by discrete-time counting processes}
\begin{document}
\maketitle

\begin{abstract}
We introduce a new class of asymmetric random walks on the one-dimensional infinite lattice.
In this walk the direction of the jumps (positive or negative) is determined by a discrete-time renewal process which is independent of the jumps. We call this discrete-time counting process the `generator process' of the walk. We refer the so defined walk to as  `asymmetric discrete-time random walk' (ADTRW). We highlight connections of the waiting-time density generating functions with Bell polynomials. 
We derive the discrete-time renewal equations governing the time-evolution of the ADTRW and analyze recurrent/transient features of simple ADTRWs (walks with unit jumps in both directions). We explore the connections of the recurrence/transience with the bias: Transient simple ADTRWs are biased and vice verse. Recurrent simple ADTRWs are either unbiased in the large time limit or `strictly unbiased' at all times with symmetric Bernoulli generator process. In this analysis we highlight the connections of bias and light-tailed/fat-tailed features of the waiting time density in the generator process. As a prototypical example with fat-tailed feature we consider the ADTRW with Sibuya distributed waiting times.

We also introduce time-changed versions: We subordinate
the ADTRW to a continuous-time renewal process which is independent from the generator process and the jumps to define the new class of 
`asymmetric continuous-time random walk' (ACTRW). This new class - apart of some special cases - is not a Montroll--Weiss continuous-time random walk (CTRW). ADTRW and ACTRW models may open large interdisciplinary fields in anomalous transport, birth-death models and others.
\\[2mm]
\noindent{\it Keywords:
\\[1mm]
Asymmetric discrete- and continuous-time random walks,
recurrence/transience, discrete-time counting process, Sibuya distribution, semi-Markov and fractional chains, Bell polynomials, light-tailed/fat-tailed waiting time distributions}
\end{abstract}
\newpage
\tableofcontents
\section{Introduction}
\label{intro}
Historically the interest in biased random walks on the integer line goes back to the classical `Gambler’s Ruin Problem' and occurred already in 1656 in a correspondence between Blaise Pascal to Pierre Fermat \cite{Pascal1656}.
A simple probabilistic version is as follows. Two gamblers $A$ and $B$ play against each other a probabilistic game multiple times.
Each game is independent of the previous ones. In a game one gambler wins a certain unit of money
whereas the other loses this amount where both gamblers start with the same amount.
If both gamblers win with the same probability $p=\frac{1}{2}$, this game is `fair', whereas
for $p\neq \frac{1}{2}$ one the gamblers has an advantage. The sequence of games stops
when one of the gamblers reaches zero money units (ruin).
The time sequence of such repeated games defines a random walk on $\mathbb{Z}$
with directed unit jumps where the first passage on zero of one of the gamblers defines `the ruin condition' (end of the game sequence). For an outline of the essential features we refer to \cite{SpitzerF1976}.
The fair case $p=\frac{1}{2}$ is an unbiased walk (as an example of a martingale \cite{Doob1953,Feller1971} and consult also \cite{Hajek1987}).

In the meantime an impressive interdisciplinary field has emerged with many variants of biased random walks with applications
in areas as varied as finance (`risk theory') and in physics sophisticated models have been developed explaining anomalous transport processes \cite{MeerschertStraka2013}. Among them asymmetric (biased) diffusion has become
a major subject with a huge amount of specialized
literature
\cite{KellyMeerschaert2019,WangBarkai2020,Angstmann-et-al2017,Hou2018,Padash2019}, just to quote a few examples.
The anomalous transport and diffusion theory is mostly based on the
continuous time random walk (CTRW) approach by Montroll and Weiss \cite{MontrollWeiss1965} where
a random walk is subordinated to an independent (continuous-time) renewal process \cite{GorenfloMainardi2008,Gorenflo2009,MainardiRobertoScalas2000,ScalasGorenfloMainardi2004,MainardiGorScalas2004}.
For fat-tailed interarrival time densities in the renewal process the resulting stochastic motion is governed by time-fractional evolution equations characterized by non-markovianity and long-time memory features. For a comprehensive overview of the wide range of models we refer to \cite{FulgerScalasGermano2008,Laskin2003,MetzlerKlafter2001,MetzlerKlafter2004,BerkowitzCortisScher2006,KutnerMasoliver2017,MeerschertSigorski2019} and the references therein. These developments have launched 
the upswing of the fractional calculus \cite{SamkoKilbasMarichev1993,OldhamSpanier1974,MainardiGorenflo1998} and generalizations \cite{Kochubei2011,Giusti-et-al-2020,Giusti2020,Luchko2021,Diethelm2020,TMM_APR_2020,TMM_APR-2020_PHYSA} (and see the references therein).

Most of the mentioned models consider continuous-time renewal processes and have profound connections
with semi-Markov chains \cite{Levy1956,Smith1955,Takacs1958,Pyke1961,Feller1964}.
In contrast, the discrete-time counterparts of semi-Markov processes, renewal processes with integer valued interarrival times
and corresponding random walk models
are relatively little touched in the literature. Essential elements of this theory have been developed only recently by Pachon, Polito and Ricciuti \cite{PachonPolitoRicciuti2021}.  
For recent pertinent physical applications in discrete-time random walks and related stochastic motions on undirected graphs we refer to our recent article \cite{MichelitschPolitoRiascos2021}.
\\[1mm]
The goal of the present paper is to introduce a new class of biased random walks where the direction of the jumps
is selected by the trials of a discrete-time counting process. We call this discrete-time counting process the `generator process' of the walk.
The approach can be extended to different directions, for instance to multidimensional biased walks.
In our model we focus on walks on $\mathbb{Z}$ and consider cases where
the asymmetry of the walk solely originates from the generator process.
\\[1mm]
The structure of our paper is as follows. In Section \ref{problem_statement} we recall
some basic mathematical features of biased walks on the integer-line to define the new class of `asymmetric discrete-time random walk' (ADTRW). We derive general expressions for the transition matrix of the `simple ADTRW' which is the ADTRW with unit jumps in both directions.
\\[1mm]
Section \ref{trial_scheme} is devoted to introduce the 'generator process' of the ADTRW. 
We define the generator process as a discrete-time counting process (renewal process with IID integer interarrival times) coming along as trial process.
To generate the ADTRW two possible outcomes of the trials ``success'' or ``fail'' determine the direction of the jumps (positive or negative).
We consider especially the long-time memory and non-markovian effects on the bias of the walk.
We introduce a scalar counterpart of the transition matrix, the `state polynomial' of the generator process which contains the complete stochastic information of the simple ADTRW (i.e. of the walk with directed unit jumps).
\\[1mm]
In Section \ref{biased-discrete-time} we highlight general connections of the waiting-time generating functions of the generator process with Bell polynomials.
\\[1mm]
Section \ref{recurrence_transience_ADTRW} is devoted to the recurrence and transience features of the ADTRW.
We analyze the connection of the memory of the waiting-time densities in the generator process with
the  recurrence/transience behavior and derive general expressions for the expected sojourn on sites 
in infinitely long simple ADTRWs.
We show for simple ADTRWs that recurrence requires light-tailed (LT) waiting time densities (short-time memory) where
the recurrent cases are
unbiased in the long-time limit. Among the recurrent simple ADTRWs it turns out that only the one with symmetric Bernoulli generator process is unbiased in a strict sense. We prove that fat-tailed (FT) waiting time densities (long-time memory) generate transient and biased simple ADTRWs.
The general approach for simple ADTRWs boils down to well known classical results in the case of the Bernoulli generator process.
\\[1mm]
As a proto-typical example with non-markovian long-memory features and FT Sibuya distributed waiting time, we introduce in Section \ref{Sibuya_walk} the `{\it Sibuya ADTRW}'. We derive the transition matrix and for the simple Sibuya ADTRW
the expected position of the walker and reconfirm the transience and bias of this walk as a consequence of the non-markovianity of the Sibuya generator process.
\\[1mm]
In Section \ref{time_changed_walk} we introduce
a time-changed version of the ADTRW: We subordinate the ADTRW to an independent continuous-time point process such as Poisson or fractional Poisson. In this way we define a new class of asymmetric
continuous-time random walks (ACTRW) which in general (apart of some special cases considered at the end of this section) are not Montroll--Weiss CTRWs. We also derive the time-evolution equations for the ACTRW which are of general fractional type.
\section{Statement of the problem and preliminary remarks}
\label{problem_statement}
In this section we recall the basic mathematical background we repeatedly use in the conception of our random walk model.
We consider a class of random walks $Y_{t \in \mathbb{N}_0} \in \mathbb{Z}$ a.s.\ on the integer line
characterized by
\beq
\label{stepn}
Y_t = \sum_{j=1}^{t} X_j ,\hspace{0.5cm} Y_0=0,
\hspace{0.5cm} X_j \in \mathbb{Z} \setminus \{0\}, \hspace{0.5cm} t \in \mathbb{N}_0
\eeq
where we allow positive and negative integer jumps ($X_j \neq 0$) taking place at integer times $t$.
We consider the initial condition that the walk starts in the origin at time $t=0$. We identify $Y_t \in \mathbb{Z}$ with the position node of a random walker at time $t$ on the infinite one-dimensional lattice.
We focus on a wider class of random walks (\ref{stepn}) which are governed by
a transition matrix of the following general type
\beq
\label{evolution_eq}
{\mathbf P}(t) = \sum_{n=0}^t \mathbb{P}(N(t)=n) ({\mathbf W}^{+})^n ({\mathbf W}^{-})^{t-n} ,\hspace{1cm} P_{i,j}(t)|_{t=0}= \delta_{i,j} ,\hspace{1cm} t \in \mathbb{N}_0
\eeq
with the elements $P_{i,j}(t)=P_{0,j-i}(t)=  \mathbb{P}(Y_t=j-i)$ ($i,j \in \mathbb{Z}$) indicating the probability that the walker is present on node $j$ at time $t$ when having started the walk on node $i$ at $t=0$ and $\delta_{i,j}$ stands for the Kronecker symbol. $\mathbf{W}^{+}$ and $\mathbf{W}^{-}$ indicate the transition matrices for positive and negative jumps, respectively.
The transition matrix (\ref{evolution_eq}) as well as $\mathbf{W}^{+}$ and $\mathbf{W}^{-}$ are `T\"oplitz' with the circulant property (defined in Eq. (\ref{circulant-def})).
In Eq. (\ref{evolution_eq}) the integer random variable
$N(t) \in \mathbb{N}_0$ (with $0\leq N(t)\leq t \in \mathbb{N}_0$) is a discrete-time counting process to be specified later (which we call the `generator process' of the walk) for the choice of the direction of the jumps. $\mathbb{P}(N(t)=n)$ indicates the probability of occurrence of $n$ positive and $t-n$ negative jumps within time interval $[0,t]$. Hence the generator process introduces (apart of some special cases) asymmetry (bias) to the walk.
We call the class of walks (\ref{stepn}), 
with transition matrix (\ref{evolution_eq}) `asymmetric discrete-time random walk' (ADTRW).
We assume that if a jump $X_j$ is positive, it is governed by the single--jump transition matrix ${\mathbf W}^{+}$ and a negative jump $X_j$ is following the single--jump transition matrix ${\mathbf W}^{-}$. The matrix ${\mathbf W}^{+}$ has the elements
\beq
\label{pos_transmt}
W^{+}_{p,q}= \Theta(q-p) {\cal W}^{+}(|q-p|)
\eeq
indicating the probability to move from $p \to q$
in one single jump.
The condition $W^{+}_{0,0}={\cal W}^{+}(0)=0$ ensures that only non-zero jumps occur. Eq. (\ref{pos_transmt}) has non-vanishing elements only in side diagonals above the main diagonal and can be seen as a right--sided discrete jump density supported
on $\ell =\{1,2,\ldots \} \in \mathbb{N}$ allowing solely strictly positive integer jumps.
In Eq. (\ref{pos_transmt}) further is introduced the `discrete Heaviside function'
defined by
\beq
\label{discrete-theta}
\Theta(r-k)= \sum_{j=-\infty}^r \delta_{j,k}  = \left\{\begin{array}{l}\ds  1 , \hspace{1cm} r-k \geq 0 \\ \\ \ds
0 , \hspace{1cm}  r-k < 0 \end{array}\right.
\eeq
where we emphasize that in our definition $\Theta(0)=1$.
Correspondingly we introduce the transition matrix for negative jumps by
\beq
\label{neq_transmt}
W^{-}_{p,q}=  \Theta(p-q) {\cal W}^{-}(|q-p|)
\eeq
with $W^{-}_{0,0}= {\cal W}^{-}(0)=0$ and is a left--sided discrete jump density.
Transition matrices
(in our convention) fulfill row-stochasticity, i.e.  
$\sum_{r=-\infty}^{\infty} W^{\pm}_{r,p} =1$ with $0\leq W^{\pm}_{r,p} \leq 1$.
We mainly consider cases where the transition matrices of Eqs. (\ref{pos_transmt}) and (\ref{neq_transmt}) have mirror symmetry $W^{-}_{0,-(q-p)}= W^{+}_{0,q-p}$, i.e. the {\it jump length} for positive and negative jumps follow the same one--sided distribution (\ref{pos_transmt}). However, bear in mind the model to be developed includes arbitrary single-jump one-sided transition matrices $\mathbf{W}^{+}$ and $\mathbf{W}^{-}$ without further symmetries.
 \\[3mm]
{\it Circulant matrices}
\\
The transition matrices and all matrices we are dealing with (including the unit matrix represented with elements
$\delta_{p,q}$, $p,q \in \mathbb{Z}$) are T\"oplitz or synonymously circulant. We call a matrix
$\mathbf{A}$ T\"oplitz if it fulfills
\beq
\label{circulant-def}
A_{p,q}= A_{p+r,q+r} ,\hspace{1cm} p,q,r \in \mathbb{Z}
\eeq
 i.e. the main diagonal and all side diagonals have identical elements, respectively.
Any circulant matrix (\ref{circulant-def}) can be represented by shift operators such that
\beq
\label{shift_p_repp}
\ds {\hat A} = \sum_{k=-\infty}^{\infty} A_{0,k}{\hat T}_{-k} , \hspace{1cm} {\mathbf A}= {\hat A}{\mathbf 1}
\eeq
where ${\mathbf 1}=(\delta_{p,q})$ denotes the unit matrix
and we introduced the spatial shift-operator ${\hat T}_r$ ($r\in \mathbb{Z}$)
which is such that ${\hat T}_{r}f(p) = f(p+r) $
($p,r \in \mathbb{Z}$) with ${\hat T}_r={\hat T}_1^r$ and we denote
with $1 = {\hat T}_0$ the zero shift.
We agree that in the notation ${\hat T}_r \delta_{p,q}=\delta_{p,q+r}=\delta_{-r,q-p}$ the shift operator acts on the right index of the Kronecker symbol and by this convention the expression ${\hat A}{\mathbf 1}={\mathbf A}$ in Eq. (\ref{shift_p_repp}) gives a well-defined circulant matrix representation.
For instance consider a single jump transition matrix generating right--sided jumps of size $+1$, namely $W^{+}_{q,p}={\hat T}_{-1}\delta_{q,p} =\delta_{q,p-1}$ and $f(p)=\delta_{i,p}$ (where the walker sits on $i$ before the jump)
to give $\sum_{r=-\infty}^{\infty}f(r)\delta_{r,p-1} =f(p-1)={\hat T}_{-1} f(p)= \delta_{i,p-1}=\delta_{i+1,p}$, i.e. the walker jumps from position $i$ to $i+1$.
The shift operator ${\hat T}_m$ ($m\in \mathbb{Z}$) is unitary and has eigenvalues $e^{i m\varphi}$ on the unit circle. Any circulant
matrix (\ref{shift_p_repp}) has canonical representation, e.g. \cite{TMM-FP-APR-Gen-Mittag-Leffler2020,RiascosMichelitschPizarro2020}
\beq
\label{canocic}
A_{m,n}=A _{0,n-m}= \frac{1}{2\pi}\int_{-\pi}^{\pi}
A(\varphi)e^{i\varphi(n-m)}{\rm d}\varphi
\eeq
with eigenvalues $A(\varphi)= \sum_{q=-\infty}^{\infty} A_{0,q}e^{-i q\varphi}$ to the (right) eigenvectors with components $\frac{e^{-i\varphi q}}{\sqrt{(2\pi)}}$ where $\varphi \in (-\pi,\pi]$.
Matrix multiplications among circulant matrices commute
 and are equivalent to discrete convolutions \cite{TMM-FP-APR-Gen-Mittag-Leffler2020} (and see Appendix in \cite{MichelitschPolitoRiascos2021} for some properties) as a consequence of commutation of shifts.
The above single--jump transition matrices can be represented by shift operators as follows
\beq
\label{shift-rep}
\begin{array}{l}
\ds {\mathbf W}^{+} = \sum_{r=1}^{\infty} {\cal W}^{+}(r) {\hat T}_{-r}{\mathbf 1}, \\[3mm]
\ds {\mathbf W}^{-} = \sum_{r=1}^{\infty} {\cal W}^{-}(r) {\hat T}_r {\mathbf 1}.
\end{array}
\eeq
Despite we focus on walks with discrete jumps, it is
rather straight-forward to extend the random walk models developed here to their continuous-space counterparts.
The single--jump transition matrices
$\mathbf{W}^{+}$ and $\mathbf{W}^{-}$ are then replaced by right-- and left-- sided continuous-space transition density kernels.

\subsection{Simple ADTRWs}
\label{simple_ADTRW}

In random walk theory an important class consists in walks where only positive and negative
jumps of unit size occur. We call this class of walks here `simple walks' or `simple ADTRWs'.
Simple walks are governed by the transition matrices ${\mathbf W}_1^{+}={\hat T}_{-1}{\mathbf 1}$ for jumps  ``$+1$''
and ${\mathbf W}_1^{-}={\hat T}_{+1}{\mathbf 1}$ for jumps ``$-1$''.
Simple random walks were extensively studied in the literature \cite{SpitzerF1976,Feller1971,PolyaG1921,RednerS,LawlerLimick2012,TMM-APR-ISTE2019}, and see the references therein. Simple walks include models for birth--death processes \cite{KarlinMcGregor1957}
with a vast field of applications in epidemiology, demography, queueing theory, finance strategies, the Gambler's Ruin Problem, and others.
In a simple ADTRW the position (\ref{stepn}) of the walker at time $t$
is given by the integer random variable
\beq
\label{simple_walk-var}
[Y_t]_{simple}= 2N(t)-t ,\hspace{1cm} Y_0=0, \hspace{1cm} t \in \mathbb{N}_0
\eeq
where
$N(t)$ jumps of $+1$ and $t-N(t)$ jumps of $-1$ are made within the time interval $[0,t]$ 
and $N(t)\in \mathbb{N}_0$ is the above-mentioned discrete-time counting process to be specified subsequently. As $0\leq N(t)\leq t$ we have that
$-t \leq [Y_t]_{simple} \leq t$.
\\[1mm]
The simple ADTRW has the drift term $-t$. It is therefore convenient to introduce a coordinate system having its origin on
a moving particle navigating with constant speed $-1$ and starting in the same position as the walker at $t=0$. The position $q(t)$ of the random walker seen by this moving particle has no drift anymore and is given by the random variable
\beq
\label{second-random}
q(t)= [Y_t]_{simple} -(-t) = 2N(t)
\eeq
corresponding to a strictly increasing walk
with positive jumps of size $2$ (almost surely).
For our convenience we introduce the transition matrix ${\mathbf Q}(t)=({\mathbf W}^{-})^{-t}{\mathbf P}(t)$ (see also Eq. (\ref{evolution_eq})) which is `seeing' the moving particle, namely 
\beq
\label{Q_t_trans_mat}
\begin{array}{cllr}
\ds {\mathbf Q}(t) & = \ds \sum_{n=0}^t\mathbb{P}(N(t)=n) ({\mathbf W}^{+})^n ({\mathbf W}^{-})^{-n} & = \ds 
\sum_{n=0}^t\mathbb{P}(N(t)=n) {\hat T}_{-2n}{\mathbf 1}    & \\ \\
\ds Q_{0,r}(t) &= \ds \sum_{n=0}^t\mathbb{P}(N(t)=n)\delta_{0,r-2n}  & = \ds  \mathbb{P}(q(t)=r) &
\end{array}
\eeq
with initial condition $P_{0,r}(t)\big|_{t=0}= Q_{0,r}(t)\big|_{t=0}= \delta_{0,r}$. The matrix
(\ref{Q_t_trans_mat})
indicates the probability that the walker at time $t$ has distance $r \in \mathbb{N}_0$ from the moving particle.
Correspondingly, the transition matrix (\ref{Q_t_trans_mat}) reduces to
\beq
\label{simple_walk_transmat}
\begin{array}{clr}
\ds Q_{0,r}(t) & = \ds \Theta(r)  \delta_{r,2\ceil{\frac{r-1}{2}}}\mathbb{P}\left(N(t)=
\frac{r}{2}\right) , \hspace{0.5cm} r \in \mathbb{Z}, \hspace{0.5cm} t \in \mathbb{N}_0 & \\ \\
\ds  P_{0,r}(t) & = \ds  Q_{0,r+t}(t) .&
 \end{array}
\eeq
$Q_{0,r}(t)$ has non-zero entries only on the sites $r \in \{0, 2,\ldots 2t-2, 2t\}$ (and $P_{0,r}(t)$ on
$r \in \{-t, -t+2,\ldots t-2, t\}$).
We introduced above the ceiling function $\ceil{a}$ indicating the
smallest integer greater or equal to $a \in \mathbb{R}$ and the Kronecker symbol picks up the terms for which $r$ is even.
Of interest is especially the probability of return to the departure site
\beq
\label{simple_walk_transmat_return}
P_{0,0}(t)  =  \delta_{t,2\ceil{\frac{t-1}{2}}} \mathbb{P}\left(N(t)=\frac{t}{2}\right)     = \left\{\begin{array}{l} \ds
\mathbb{P}\left(N(t)=\frac{t}{2}\right) , \hspace{0.5cm}  \hspace{0.5cm} t \hspace{0.25cm} {\rm even} \\ \\
\ds 0 ,\hspace{0.5cm} t \hspace{0.25cm}  {\rm odd}. \end{array}\right.
\eeq
In the subsequent section we specify the trial process selecting the direction of the jumps in Eq. (\ref{evolution_eq}) and recall the notion of `discrete-time counting process' which comes along as `renewal trial process' with connections to discrete-time semi-Markov chains. For essential elements of the theory consult \cite{PachonPolitoRicciuti2021}.

\section{Discrete-time renewal processes and trial schemes}
\label{trial_scheme}
Here we introduce a trial process which defines the directions of the jumps in the ADTRW. 
Consider a sequence of trials
where each trial has two possible outcomes, ``success'' or ``fail''. For our convenience we introduce random variables $Z_r\in \{0,1\}$, a.s., representing these possible outcomes, namely $\{Z_r=0\}$ for a fail and
$\{Z_r=1\}$ for a success at trial $r \in \mathbb{N}$. Furthermore, let $\mathbb{P}(Z_r=1|Z_{r-1}=0)=\alpha_r \in [0,1]$, $r \in \mathbb{N}$, $\mathbb{P}(Z_r=1|Z_{r-1}=1)=\alpha_1$, be the conditional probability of success in the $r$-th trial conditional to the filtration $\mathcal{F}_{r-1}$ (up to trial $r-1$) to which the trial process is adapted. Performing a sequence of $k$ trials gives $2^k$ possible outcomes. Each outcome $(z_1,z_2,\ldots, z_k)$ occurs with probability $P_k(z_1,z_2,\ldots,z_k)$.
The probability of a certain outcome in a sequence of $k$ trials has the structure
\beq
\label{independence}
P_k(z_1,z_2,\ldots,z_k)= p_1(z_1) p_2(z_2)\ldots p_k(z_k),
\eeq
where the $p_r(z_r)$, $r \in \{1,\dots,k\}$ are in general conditional probabilities.
As an example, the probability that a trial sequence
``(fail, fail, fail, success, fail, success)'' occurs in a sequence of $6$ trials then is with above adaption rule $P_6(0,0,0,1,0,1)=(1-\alpha_1)(1-\alpha_2)(1-\alpha_3)\alpha_4(1-\alpha_1)\alpha_2$. If the $\alpha_k$ is non constant the process has a memory (i.e.\ it has non geometric waiting times).
A sequence of $k$ trials where the first success occurs at trial $k$ has then the probability
\beq
\label{first-sucess}
 \psi_k = P_k(0,\ldots 0,1) = \alpha_k(1-\alpha_{k-1})\ldots(1-\alpha_1) ,\hspace{0.5cm} \alpha_j \in [0,1], \: k \in \mathbb{N}
\eeq
with the `survival probability' (probability that in $k-1$ trials all outcomes are fails) $P_{k-1}(0,\ldots,0)= S_{k-1} = \prod_{j=1}^{k-1}(1-\alpha_j)$. On the other hand we observe that any discrete density $\psi_k$ ($k\in \mathbb{N}$) can be represented as in Eq. (\ref{first-sucess}) within such a trial scheme with $S_{k-1}=\sum_{r=k}^{\infty}\psi_r$ and $S_k=S_{k-1}-\psi_k= S_{k-1}(1-\alpha_k)$. Notice that the survival probability\footnote{We do not consider here cases where $S_{\infty} =1-\sum_{k=1}^{\infty} \psi_k >0$ where a finite survival probability exists, i.e. a finite probability that in infinitely many trials never a success occurs.} $S_k \to S_{\infty}= 0$ as $k\to \infty$. With the initial condition $S_0=1$ it follows then $\sum_{k=1}^{\infty}\psi_k=1$.

A pertinent example of the form (\ref{first-sucess}) with memory is the Sibuya distribution (also referred to as Sibuya($\alpha$) which has $\alpha_k=\beta/k$, $\beta \in (0,1)$,
thus the probability of the first success at trial $k$ is
\beq
\label{Sibuyaalpha}
\begin{array}{l}
\ds \psi_{\beta}(k)  = \frac{\beta}{k}\left(1-\frac{\beta}{k-1}\right )\ldots (1-\beta) = (-1)^{k-1}\binom{\beta}{k} ,\hspace{0.5cm} \beta \in (0,1) ,\hspace{0.5cm} k\in \mathbb{N}
\end{array}
\eeq
which is fat-tailed, i.e.\ for $k \to \infty$ we have a heavy power-law tail $\psi_{\beta}(k) =\frac{\beta}{k}S_{\beta,k-1} \sim \frac{\beta}{\Gamma(1-\beta)} k^{-1-\beta}$. Hence, the Sibuya survival
probability tends to zero as a power-law
$S_{\beta}(k) =  (-1)^{k}\binom{\beta -1}{k}  \sim \frac{k^{-\beta}}{\Gamma(1-\beta)}$
reflecting the long-memory and non-markovian feature of the Sibuya trial process. We will come back to the Sibuya distribution later on.
\subsection{Discrete-time renewal process - generator process of the ADTRW}
\label{generator_process}
In order to establish the connection with discrete-time renewal processes with above mentioned adaption rule we consider the conditional probabilities
\beq
\label{ren_def}
\begin{array}{clr}
\ds \mathbb{P}(Z_k=1|Z_{k-1}=0) & = \ds \alpha_k , \hspace{1cm}&  \alpha_k \in [0,1] 
 \\ \\
\ds \mathbb{P}(Z_k=1|Z_{k-1}=1) & = \ds \alpha_1, & k \geq 2
\end{array}
\eeq
and probability of success $\mathbb{P}(Z_1=1)= \alpha_1$ in the first trial.
Once in a trial sequence a success occurs (i.e. $Z_r=1$) for the first time, say at trial $r$, then the probability for success in the subsequent trial is reset to $\alpha_1$ as in the first trial and the process starts anew. In the renewal picture the trial number $r$ of first success can be seen as the integer renewal time (arrival time of event `success') in a discrete-time renewal process.
\\[1mm]
With these remarks we introduce the discrete-time counting process $N(t)$ which counts the successes among $t$ trials (success $=$ `arrival' or `event' in the counting process) as
\beq
\label{counting_rocess}
N(t)= \sum_{j=1}^t Z_j  , \hspace{1cm} N(0)=0, \hspace{1cm} Z_j \in \{0,1\} \text{ a.s} .
\eeq
We call the so defined discrete-time counting process with conditional probabilities (\ref{ren_def}) `generator process' of the ADTRW.
The number of trials between two successes then can be seen as IID interarrival times $\Delta t_j$ which can take positive integer values $\Delta t_j \in \{1,2,\dots \}$. Further, to define the jump directions in the ADTRW, we associate with each `success' a positive jump 
and with each `fail' a negative jump. Hence, the integer variable $N(t) \in \mathbb{N}_0$ counts the number of positive jumps and $t-N(t)$ the number of negative jumps within the time interval $[0,t]$ and clearly we have $0\leq N(t) \leq t$.

Then we define a discrete-time renewal process fulfilling property (\ref{ren_def}) where the IID interarrival times follow the waiting-time density ($\Delta t_j \to t$)
\beq
\label{define_renewal}
 \mathbb{P}(Z_1=0,\ldots,Z_{t-1}=0,Z_t=1) = \psi(\alpha_1,\ldots,\alpha_t) = \alpha_t \prod_{j=1}^{t-1}(1-\alpha_j) ,\hspace{1cm} t\in \{1,2\ldots\}.
\eeq
We use from now on the synonymous notations
$\psi_t=\psi(t)=\psi(\alpha_1,\ldots,\alpha_t)$
for Eq.\,(\ref{define_renewal}). 
Of utmost importance are the `state probabilities', i.e. the probabilities that in $t$ trials $n$ successes occur ($n$ arrivals up to time $t$). The state probabilities are defined as
\beq
\label{ste_prob}
\Phi^{(n)}(\alpha_1,\ldots,\alpha_t) = \mathbb{P}(N(t)=n) = \mathbb{P}\left(\sum_{j=1}^t Z_j =n\right)  ,\hspace{0.5cm} 
n, t \in \{0,1,2,\ldots\} \eeq
In general, for non-constant $\alpha_t$ the complete history of the outcomes of $t$ trials is considered, and the generator process has a memory and is non-Markovian. 
The probability of no success ($t$ successive fails) in $t$ trials
writes
\beq
\label{survival_prob}
\begin{array}{clr}
\ds \mathbb{P}(N(t)=0) &= S(\alpha_1,\ldots,\alpha_t) = \prod_{j=1}^t(1-\alpha_t)
, & \ds t \in \mathbb{N}\\ \\
  \mathbb{P}(N(t)=0)\big|_{t=0} &  =\ds  1 , &
\end{array}
\eeq
i.e.\ the survival probability which we have previously introduced. We have the initial condition
 $\mathbb{P}(N(t)=n)|_{t=0}= \delta_{n,0}$ and we point out a further feature of  discrete-time counting processes, namely 
 $\mathbb{P}(N(t)=n)$ is non-null only for $n \in [0,t]$,
 reflecting that $0\leq N(t) \leq t$.
One further observes the normalization condition
\beq
\label{normalization_state}
  \mathbb{P}(N(t) \leq t) = \sum_{n=0}^t \mathbb{P}(N(t)=n) = 1,
\eeq
of the state distribution (\ref{ste_prob}). In other words Eq. (\ref{normalization_state}) covers all $2^t$ possible paths in the branching tree of the $t$ trials.
In order to connect the trial process with the above biased random walk (see Eq. (\ref{evolution_eq}))
we define the transition matrix for the jump taking place at instant $t \in \mathbb{N}$ as
\beq
\label{transmat_time_t}
{\mathbf W}_t(Z_t) = Z_t{\mathbf W}^{+}+ (1-Z_t){\mathbf W}^{-} ,\hspace{0.5cm} Z_t \in \{0,1\},\hspace{0.5cm} t \in \mathbb{N}.
\eeq
For $Z_t=1$ (`success') the walker makes a positive jump following ${\mathbf W}^{+}$ and a negative jump following
${\mathbf W}^{-}$ otherwise.
Then we have that
\beq
\label{write_4_trans_mat}
\begin{array}{clr}
\ds {\mathbf P}(t) & = \ds   \mathbb{E}\prod_{j=1}^t [Z_j{\mathbf W}^{+} +(1-Z_j){\mathbf W}^{-}] & 
\\ \\
 & = \ds  \mathbb{E} \left[({\mathbf W}^{+})^{N(t)} ({\mathbf W}^{-})^{t-N(t)}\right] & 
\\ \\
 & = \ds  \sum_{n=0}^t\mathbb{P}(N(t)=n) ({\mathbf W}^{+})^n ({\mathbf W}^{-})^{t-n} &
 \end{array}
\eeq
which is the initially claimed ADTRW transition matrix (\ref{evolution_eq}).
\\[1mm]
For our convenience we make extensively use of generating functions.
Let $\psi(t)$ be a discrete-time density of Eq. (\ref{define_renewal}). Then we
introduce its generating function by\footnote{We indicate generating functions of densities $f(t)$ by
${\bar f}(u)$.}
\beq
\label{genfu}
{\bar \psi}(u)  =\sum_{t=1}^{\infty}\psi(t) u^t ,\hspace{1cm} |u| \leq 1
\eeq
with $\psi(t)= \frac{1}{t!}\frac{d^t}{du^t}{\bar \psi}(u)|_{u=0}$ and where ${\bar \psi}(u)|_{u=1}=1$ reflects the normalization of density (\ref{define_renewal}).
Further, we impose
in generating function (\ref{genfu}) the initial condition $\psi(t)|_{t=0}=0$ ensuring that the minimum waiting-time between two successes is $\Delta t=1$.
Then we introduce the convolution of two discrete distributions $g(t), h(t)$ supported on $\mathbb{N}_0$ by
\beq
\label{dconvolutions}
[g \star h] (t) =: \sum_{n=0}^t g(n)h(t-n)
\eeq
with generating function $\sum_{t=0}^{\infty} u^t [g \star h] (t) = {\bar g}(u){\bar h}(u)$.
We
denote convolution powers as $$[\underbrace{g\star \ldots \star g}_{n \, \, times}](t) = [g \star]^n(t),$$ with 
generating functions $({\bar g}(u))^n$ (where especially $[g \star]^1(t)=g(t)$ and $[g \star]^0(t)=\delta_{t,0}$).
For an outline of the connections between generating functions, discrete-time convolutions and related shift-operators we refer to our recent article \cite{MichelitschPolitoRiascos2021}.
By simple conditioning arguments one obtains for the state probabilities (\ref{ste_prob})
\beq
\label{state_probs}
\mathbb{P}(N(t)=n) =\Phi^{(n)}(t) = \Phi^{(0)}(t) \star [\psi(t)\star]^n =
\left[S(\alpha_1,\ldots,\alpha_t) \star [\psi(\alpha_1,\ldots,\alpha_t) \star]^n\right](t) ,\hspace{0.5cm} n,t \in \mathbb{N}_0.
\eeq
The state probabilities (\ref{state_probs}) are non-zero for $0 \leq n \leq t$ simply telling us that the number of successes in $t$ trials is within
$0\leq N(t)\leq t$.
Especially convenient is to employ the generating function of the state probabilities (see \cite{PachonPolitoRicciuti2021,MichelitschPolitoRiascos2021} for detailed derivations) 
\beq
\label{gen_state-function}
{\bar \Phi}^{(n)}(u) = \sum_{t=0}^{\infty}\mathbb{P}(N(t)=n) u^t = \sum_{t=n}^{\infty}\Phi^{n}(t) u^t
=\frac{1-{\bar \psi}(u)}{1-u} ({\bar \psi}(u))^n ,\hspace{0.5cm} n \in \mathbb{N}_0
\eeq
with $\mathbb{P}(N(t)=n)= \frac{1}{t!}\frac{d^t}{du^t}{\bar \Phi}^{(n)}(u)|_{u=0}$ where ${\bar \psi}(u)=\alpha_1 u+o(u) = u  \sum_{t=1}^{\infty}\psi(t)u^{t-1}$
having lowest order $u$ reflecting ${\bar \psi}(u)|_{u=0} = \psi(t)|_{t=0}=0$. Hence
${\bar \Phi}^{(n)}(u)= \alpha_1^nu^n +o(u^n)$ thus confirms
$\mathbb{P}(N(t)=n)=0$ for $t<n$.
We further observe that $\mathbb{P}(N(t)=n)|_{t=n}= \alpha_1^n$ which is the situation when in $t$ trials all outcomes are successes with $N(t)=t$.
It is convenient to introduce the polynomial of degree $t$ (generating function of the state probabilities)
\beq
\label{PZ}
{\cal P}(v,t)= \mathbb{E} v^{N(t)} =\sum_{n=0}^{\infty} v^n \mathbb{P}(N(t)=n) =\sum_{n=0}^t v^n \mathbb{P}(N(t)=n) ,\hspace{1cm} t \in \mathbb{N}_0
\eeq
where the series stops at $n=t$ as $\mathbb{P}(N(t)=n)=0$ for $n>t$. Thus this generating function is a polynomial of order $t$.
We call ${\cal P}(v,t)$ the `state polynomial' of the generator process.
Useful is also a rescaled version (see Eq. (\ref{write_4_trans_mat}))
\beq
\label{similar}
\Lambda(a,b,t) =:
%\mathbb{E}[a Z +b(1-Z)] =
\mathbb{E} \,[ a^{N(t)} b^{t-N(t)} ] = \sum_{n=0}^t a^nb^{t-n} \mathbb{P}(N(t)=n) =  
b^t \, {\cal P}\left(\frac{a}{b},t\right).
\eeq
We observe that ${\cal P}(v,t)\big|_{v=1}=\Lambda(1,1,t)=1$ as a consequence of the normalization
(\ref{normalization_state}) and we have $\Lambda(v,1,t)={\cal P}(v,t)$ with the scaling property $\Lambda(\lambda a,\lambda b,t)=\lambda^t\Lambda(a,b,t)$.
The state polynomials (\ref{PZ}) and (\ref{similar})
in the ADTRW model come along as matrix functions defining the transition matrix (\ref{write_4_trans_mat}): ${\mathbf P}(t)= \Lambda({\mathbf W}^{+},{\mathbf W}^{-},t)$.
The state polynomial contains the complete stochastic information of the simple ADTRW.
So for instance the expected position of the walker in a simple ADTRW at time $t$ is obtained as
\beq
\label{simple_walk_bias}
\begin{array}{clr}
\ds \mathbb{E}[Y_t ]_{simple} &= \ds \mathbb{E}[N(t) - (t-N(t)) ]  = \left(\frac{\partial }{\partial a}- \frac{\partial }{\partial b} \right) \Lambda(a,b,t)\Big|_{a=b=1}   & \\[3mm]
 &= \ds 2 \frac{\partial}{\partial v}{\cal P}(v,t)\Big|_{v=1} -t,&
 \end{array}
\eeq
containing the drift term $-t$ which is removed
from the coordinate system of the moving particle:
$\mathbb{E}[q(t)] = 2 \frac{\partial}{\partial v}{\cal P}(v,t)\big|_{v=1} = 2 \mathbb{E}[N(t)]$. This corresponds to a strictly increasing walk with
jumps of size $2$, see Eq. (\ref{second-random}). Clearly, the expected position of the walker in a simple ADTRW is bounded $-t \leq \mathbb{E} [Y_t]_{simple} \leq t$. 
We point out that we employ the synonymous notation 
$\mathbb{E}[A] = \mathbb{E} A$ for expectation values of random variables $A$ where we sometimes omit the braces $[..]$.
Then, it is convenient to introduce the generating function of the state polynomial
\beq
\label{gen_state_pol}
\begin{array}{clr}
\ds {\bar {\cal  P}}(v,u) =: \sum_{t=0}^{\infty} u^t {\cal P}(v,t) & = \ds
\sum_{t=0}^{\infty} u^t \sum_{n=0}^{\infty} v^n \mathbb{P}(N(t)=n) ,\hspace{0.5cm} |u| < 1 ,
\hspace{0.5cm} |v| \leq 1& \\ \\
\ds & = \ds \sum_{n=0}^{\infty}v^n {\bar \Phi}^{(n)}(u) = 
\frac{1}{(1-u)}\frac{1-{\bar \psi}(u)}{1-v{\bar \psi}(u)} = \frac{{\bar \Phi}^{(0)}(u)}{1-v{\bar \psi}(u)}  &
\end{array}
\eeq
which is related to the generating function of expression (\ref{similar}), yielding
\beq
\label{similar_gen}
\begin{array}{clr}
\ds {\bar \Lambda}(a,b,u)& = \ds  {\bar {\cal  P}}\left(\frac{a}{b},ub\right)  =  \sum_{t=0}^{\infty} u^t \, \sum_{n=0}^t\mathbb{P}(N(t)=n) a^nb^{t-n} =
\frac{1-{\bar \psi}(bu)}{1-bu} \sum_{n=0}^{\infty} \frac{a^n}{b^n}({\bar \psi}(bu))^n &\\ \\
\ds  & = \ds \frac{b[1-{\bar \psi}(bu)]}{(1-bu)[b-a{\bar \psi}(ub)]}
, \hspace{1cm} |u| < 1 , \,\, |a|,|b| \leq 1,\qquad (b \neq 0) & 
\end{array}
\eeq
converging at least for $|u|< 1$. See also Appendix \ref{Appendix_A} for some pertinent limiting cases.
We have ${\bar \Lambda}(1,1,u)={\bar {\cal P}}(1,u) =\frac{1}{1-u}$ reflecting the normalization (\ref{normalization_state}) and the initial condition
\beq
\label{initial_condition}
{\bar \Lambda}(a,b,u)|_{u=0}=
\Lambda(a,b,t)|_{t=0} = {\bar {\cal P}}(v,u)|_{u=0}=
{\cal P}(v,t)|_{t=0} =1,
\eeq
as $\mathbb{P}(N(t)=n)\big|_{t=0} =\delta_{n0}$.
The ADTRW transition matrix (\ref{evolution_eq}) is then given by the matrix function
\beq
\label{evolution}
{\mathbf P}(t) = \Lambda({\mathbf W}^{+},{\mathbf W}^{-},t) = \frac{1}{t!}\frac{d^t}{du^t} \left(
\frac{{\mathbf W}^{-}[{\mathbf 1}-{\bar \psi}({\mathbf W}^{-} u)]}{(1-{\mathbf W}^{-}u)[{\mathbf W}^{-}-{\mathbf W}^{+}{\bar \psi}({\mathbf W}^{-} u)]}\right)\bigg|_{u=0}
\eeq
and fulfills the initial condition ${\mathbf P}(t)\big|_{t=0}= {\mathbf 1}$ as a consequence of Eq. (\ref{initial_condition}).
The convergence of the matrix generating function ${\bar \Lambda}({\mathbf W}^{+},{\mathbf W}^{-},u) $ is ensured by the
fact that the eigenvalues $W^{\pm}(\varphi)$
of the transition matrices ${\mathbf W}^{\pm}$
fulfill $|W^{\pm}(\varphi)| \leq 1$ \cite{TMM-FP-APR-Gen-Mittag-Leffler2020,RiascosMichelitschPizarro2020}.
In general we have that ${\bar \Lambda}(a,b,u)\neq 
{\bar \Lambda}(b,a,u)$ where the absence of the exchange symmetry is telling us that the ADTRW is in general biased even if $W^{+}_{0,p-q}= W^{-}_{0,-(p-q)}$ have mirror symmetry. As a consequence of the asymmetry of the walk the eigenvalues of the transition matrix $\Lambda(e^{-i\varphi},e^{i\varphi},t)$ are complex and the transition matrix (\ref{evolution}) is generally not symmetric. 
This holds true with the exception of the class of `strictly unbiased walks' where the transition matrices are self-adjoint (symmetric) with real eigenvalues being even functions of $\varphi$. We will prove in
Section \ref{recurrence_transience_ADTRW}
that a simple ADTRW is strictly unbiased only if its generator process is the symmetric Bernoulli process.
\\[1mm]
We observe that the state polynomial $\Lambda(a,b,t)$ fulfills the following renewal equation (see Appendix
\ref{Appendix_B})
\beq
\label{renewal_eq_Lam}
\begin{array}{clr}
\ds \Lambda(a,b,t) & = \ds  b^t \Phi^{(0)}(t)  + \sum_{r=1}^t a b^{r-1} \psi(r) \Lambda(a,b,t-r), \hspace{0.5cm} t \in \mathbb{N}, & \\ \\
\ds \Lambda(a,b,t)\big|_{t=0} & = \ds 1
\end{array}
\eeq
containing the survival probability $\mathbb{P}(N(t)=0)= \Phi^{(0)}(t)=S(\alpha_1,\ldots,\alpha_t)= \prod_{\ell=1}^t(1-\alpha_{\ell})$ and the waiting time density
$\psi(t)=\alpha_t S(\alpha_1,\ldots,\alpha_{t-1})$ (see Eqs. (\ref{define_renewal}), (\ref{survival_prob})).
The right-hand side of the renewal equation contains the history of the process $\{\Lambda(a,b,r)\}$ ($0 \leq r \leq t-1$) reflecting the memory and the non-markovian nature of the ADTRW.
The renewal equation is especially useful for numerical evaluations to successively compute $\Lambda(a,b,t)$ from all its previous values and the waiting time density $\psi(r)$ ($r\in \{1,\ldots,t\}$) of the generator process.
For instance, for $t=1$ we have
$\Lambda(a,b,1)= b(1-\alpha_1)+a\alpha_1$, and so forth.
The renewal equation
for the state polynomial is contained in Eq. (\ref{renewal_eq_Lam}) by accounting for ${\cal P}(v,t)=\Lambda(v,1,t)$. 
By simply replacing $a \to {\mathbf W}^{+}$, $b \to {\mathbf W}^{-}$ in Eq. (\ref{renewal_eq_Lam}) gives the time-evolution 
equation with memory which governs the transition matrix (\ref{evolution}).
For the simple walk
(i.e. ${\hat W}^{+}={\hat T}_{-1}$, ${\hat W}^{-}={\hat T}_{+1}$) we get
\beq\label{simple_walk_renewal}
\begin{array}{clr}
\ds P_{i,j}(t) & = \ds \Phi^{(0)}(t) \, \delta_{i,j+t}  + \sum_{r=1}^t\psi(r)P_{i,j+r-2}(t-r) , \hspace{0.5cm} t \in \mathbb{N} & \\ \\
 \ds P_{i,j}(t)\big|_{t=0} & = \ds \delta_{ij}
 \end{array}
\eeq
being solved by the transition matrix of the simple ADTRW (see Eq. (\ref{simple_walk_transmat})).
Then we can rewrite the renewal equation as a master equation with memory as
\beq
\label{master_eq_memory}
P_{ij}(t)-P_{ij}(t-1) =  \Phi^{(0)}(t) \, \delta_{i,j+t} + \sum_{r=1}^{\infty}\psi(r)[P_{i,j+r-2}(t-r)-P_{ij}(t-1)], \hspace{0.3cm} t \in \mathbb{N}
\eeq
with initial condition $P_{ij}(t)\big|_{t=0}=\delta_{ij}$ and recall causality, i.e.
$P_{ij}(t)=0$ for $t< 0$.  
Consult also Appendix \ref{Appendix_B} for some operator representations. It appears instructive to consider the following two limiting cases which correspond to strictly increasing and decreasing walks, respectively.
\\[3mm]
{{\bf (i)} {\it Limit of short waiting times: `Markovian limit'}
\\[2mm]
A markovian limit is obtained for the `trivial case' when each trial almost surely is a success with $\psi(t)= \delta_{t,1}$ 
($\alpha_1=1$), i.e.\ the waiting time density has the shortest possible tail of one time unit.
The survival probability then is  $\Phi^{(0)}(t)=\delta_{t,0}$.
Then renewal equation (\ref{renewal_eq_Lam}) boils down to the memoryless recursion
\beq
\label{renewal_eq_Lam_markovian}
  \Lambda_{\alpha_1=1}(a,b,t) = \delta_{t,0}+  a \Lambda(a,b,t-1)
\eeq
which has, with initial condition $\Lambda(a,b,t)\big|_{t=0}=1 $, the simple solution $\Lambda_{\alpha_1=1}(a,b,t)=a^t$ independent of $b$ and therefore coincides with the limiting case ${\bar \Lambda}(a,0,u)= \frac{1}{1-au}$ for $\alpha_1=1$ and ${\bar \psi}(u)=u$ (see Eq. (\ref{limit_b_0}), Appendix \ref{Appendix_A}).
In this limit the evolution equation takes the form
\beq
\label{markovian_limit}
\begin{array}{clr}
\ds P_{ij}(t) & = \ds \delta_{ij} \delta_{t,0} + P_{i,j-1}(t-1) & \\ \\
\ds {\mathbf P}(t) & = \ds {\mathbf 1}\delta_{t,0} +{\mathbf W}^{+}
{\mathbf P}(t-1) , & \ds \hspace{1cm} W^{+}_{ij} = {\hat T}_{-1}\delta_{ij}=\delta_{i,j-1}
\end{array}
\eeq
being solved by ${\mathbf P}(t)= [{\mathbf W}^{+}]^t={\hat T}_{-t}{\mathbf 1}$ with entries
\beq
\label{markov_limit_sol}
P_{ij}(t)={\hat T}_{-t} \delta_{i,j}=\delta_{i,j-t}
\eeq
corresponding to a strictly increasing walk with unit jumps $+1$ (almost surely).
\\[3mm]
{\bf (ii)} {\it Limit of long waiting times: `Frozen limit'}
\\[2mm]
Another limiting case of interest is when the waiting time density is concentrated at
infinity, i.e.\ 
the probability for a success in the generator process becomes smaller and smaller (though not zero) and the waiting time density then has an extremely long tail. We then have that $\alpha_t \le \epsilon \to 0$ for each finite $t$. Examples for this limit include the geometric waiting-time density
$\psi_B(t)=pq^{t-1}$ for $p = \epsilon \to 0$, or in the Sibuya density (\ref{Sibuyaalpha}) this limit is obtained for $\beta \to 0$ (considered subsequently). 
In the Bernoulli case the survival probability $\Phi^{(0)}(t) = (1-p)^t$ remains close to one during a
`very long' time interval, i.e.\ at time scales $0\leq t < 1/p^{\delta}$ with $1/p^{\delta} \to \infty$ for any $\delta \in (0,1)$ and
$1 \ll 1/p^{\delta} \ll 1/p$. Indeed, for the Bernoulli ADTRW
we have for $p \to 0$ the state polynomial (\ref{Lambda_bern}) $\lim_{p\to 0+}\Lambda_B(a,b,t)=\lim_{p\to 0+}(ap+qb)^t = b^t$. This limit is therefore connected with the limit
 $\Lambda(a,b,t)$ when $a\to 0$ having generating function (\ref{limit_a_zero}), see Appendix \ref{Appendix_A}.
The effect is that the
survival probability $\Phi^{0}(t) =(1-\alpha_1)\ldots (1-\alpha_t) \to 1-$ (for $t$ finite) remains for a very long time close to its initial value one thus the walker remains a long time `frozen' in its `ground state' $N(t)=0$ (though not `forever' as eventually $\Phi^{0}(t) \to 0$ for $t \gg 1/\epsilon$).
The renewal equation for the frozen limit
becomes
\beq
\label{frozen}
 \Lambda(a,b,t) =  b^t \Phi^{(0)}(t) \sim b^t ,\hspace{1cm} 0\leq  t  < \frac{1}{\epsilon^{\delta}} \ll \frac{1}{\epsilon} , \hspace{1cm} \delta \in (0,1)
 \eeq
and it is independent of $a$, so that in such a walk negative jumps strongly dominate. The renewal equation (\ref{simple_walk_renewal}) then writes
\beq
\label{frozen-trans}
P_{ij}(t) = {\bar T}_t \delta_{i,j} \Phi^{(0)}(t) = \delta_{i,j+t}, \hspace{1cm} \Phi^{(0)}(t) 
= 1-  ,\hspace{1cm} 0\leq  t  < \frac{1}{\epsilon^{\delta}} \ll \frac{1}{\epsilon} , \hspace{1cm} \delta \in (0,1) .
\eeq
In the coordinate system of the moving particle the walker hence does not move for a long time which is reflected by
$Q_{i,j}(t)={\hat T}_{-t}P_{i,j}(t)=\delta_{i,j}$.
\section{Connections with Bell polynomials}
\label{biased-discrete-time}
Here we point out an interesting connection with a certain class of Bell polynomials \cite{Bell1928,Bell1934}.
Consider first the generating function representation of convolution powers of the waiting-time density:
\beq
\label{n_trials}
\begin{array}{cclr}
\ds [\psi *]^n(r)&=&  \ds \frac{1}{r!} \frac{d^r}{du^r}({\bar \psi}(u))^n\big|_{u=0} &  \\ \\
  &=& \ds \frac{1}{r!} \frac{d^r}{du^r}\left(\psi_1u+\ldots +\psi_{r-n+1}u^{r-n+1}\right)^n\big|_{u=0}  \\ \\
 &=& \ds \frac{1}{r!} \sum_{n_1+n_2+\ldots n_{r-n+1} = n} \frac{n!}{n_1!\ldots n_{r-n+1}!}\psi_1^{n_1}\psi_2^{n_2} \ldots
  \psi_{r-n+1}^{n_{r-n+1}} \frac{d^r}{du^r} u^{n_1+2n_2+\ldots (r-n+1)n_{r-n+1}}\big|_{u=0}  & \\ \\
 &=& \ds \sum_{n_1+n_2+\ldots n_{r-n+1} = n} \frac{n!}{n_1!\ldots n_{r-n+1}!}\psi_1^{n_1}\psi_2^{n_2} \ldots
  \psi_{r-n+1}^{n_{r-n+1}} \, \delta_{r, n_1+2n_2+\ldots (r-n+1)n_{r-n+1}} & \\ \\
&=& \ds {\cal B}_{r,n}(\psi_1,\psi_2,\ldots, \psi_{r-n+1}) ,\hspace{1cm} 1 \leq n \leq r \in \mathbb{N}. 
\end{array}
\eeq
We have $[\psi *]^n(r)=0$ for $n>r$ and therefore
\beq
\label{extend_Bell}
{\cal B}_{r,n} = 0 , \hspace{0.5cm}  n > r.
\eeq
In Eq. (\ref{n_trials}), for the expression of ${\bar \psi}(u)$ we omit the terms
with orders $t > r-n+1$ as they give zero contribution. We further have
${\cal B}_{r,1}= \psi_r$.
The trivial cases $n=0$ and $r=0$ (considering
$[\psi *]^0(r) =\delta_{r,0}$, $r\in \mathbb{N}_0$), yield
\beq
\label{extend_def}
\begin{array}{l}
\ds  {\cal B}_{r,0} = \delta_{r,0}, \hspace{1cm}    r\in \mathbb{N}_0 \\ \\
 \ds  {\cal B}_{0,n} = ({\bar \psi}(u))^n\bigg|_{u=0} = \delta_{0,n}.
 \end{array}
 \eeq
The quantities ${\cal B}_{r,n}(\psi_1,\psi_2,\ldots, \psi_{r-n+1})$ are referred to as the `incomplete ordinary Bell polynomials' \cite{Bell1928,Bell1934}.
The Kronecker symbol $\delta_{r, n_1+2n_2+\ldots (r-n+1)n_{r-n+1}}$ indicates that the only
terms which contribute are those for which $\sum_{k=1}^{r-n+1} k n_k = r$
and $n_k$ ($0\leq n_k \leq n$) are non-negative integers such that in above multinomial $\sum_{k=1}^{r-n+1}n_k=n$.
This summation covers all possible partitions of the
integer $r$ into $n$ members where each member is of integer size $k=1,2,\ldots, r-n+1 \in \mathbb{N}$. The member $k$ occurs with
multiplicity $n_k$ where for $n>r$ no such partition exists thus property (\ref{extend_Bell}) holds true.
For instance, when $r=n$ there is only one partition into $n$ members, namely each member of size $k=1$ with multiplicity $n_1=r$. On the other hand, for $n=1$ there is only one partition (i.e. $n_k=1$) of size
$k=r$.
It follows hence from Eq. (\ref{n_trials}) that
\beq
\label{inverse}
({\bar \psi}(u))^n = \sum_{t=n}^{\infty} u^t {\cal B}_{t,n}(\psi_1,\psi_2,\ldots, \psi_{t-n+1}).
\eeq
We point out that the incomplete ordinary Bell polynomials ${\cal B}_{r,n}$ and the `incomplete exponential Bell polynomials' $B_{r,n}^{exp}$ are related by \cite{Bell1934}
\beq
\label{Bexp}
\ds B_{r,n}^{exp}(x_1,\ldots x_{r-n+1})  = \ds
\frac{r!}{n!} {\cal B}_{r,n}\left(\frac{x_1}{1!},\frac{x_2}{2!},\ldots ,
\frac{x_{r-n+1}}{(r-n+1)!}\right).
\eeq
Indeed the exponential Bell polynomials come into play in the remarkable Fa\`a di Bruno's formula \cite{Johnson2002} emerging in a composition of functions from the chain rule $\frac{d^t}{dx^t}f(\tau g(x))\big|_{x=0}$.
For an outline of this beautiful theory and its interpretations in combinatorics
we refer the interested reader to
\cite{Cvijovic2011,Kruchinin2011} and the references therein.

For a fixed $t \in \mathbb{N}$, by means of the set of ordinary incomplete Bell polynomials $\{{\cal B}_{t,n}\}$ ($1\leq n\leq t$) we can generate the complete ordinary Bell polynomial as
\beq
\label{Bell_polynomial}
\begin{array}{l}
{\cal B}_{t}(\psi_1,\psi_2,\ldots, \psi_t; v)  \\ \\
=\ds
\sum_{n=1}^t v^n {\cal B}_{t,n}(\psi_1,\psi_2,\ldots, \psi_{t-n+1})
  = \frac{1}{t!}\frac{d^t}{du^t} \left(\frac{v{\bar \psi}(u)}{1-v{\bar \psi}(u)}\right)\bigg|_{u=0} ,
\hspace{0.3cm} t \in \{1,2,\ldots\}
\end{array}
\eeq
and with ${\cal B}_{0,0}=1$ we
have ${\cal B}_0=1$.
Then, by accounting for Eq. (\ref{gen_state-function}) the state polynomial writes
\beq
\label{state_polynmial_bell}
{\cal P}(v,t) =  \sum_{r=0}^t  \Phi^{(0)}(t-r)
{\cal B}_{r}(\psi_1,\psi_2,\ldots, \psi_r; v),
\eeq
where in this convolution ${\cal B}_0=1$.
Using the representation (\ref{Bell_polynomial}) we can write for the expected number of arrivals
within $[0,t]$
\beq
\label{mean_arrivals}
\ds \mathbb{E} [N(t)]  =  \frac{\partial}{\partial v} {\cal P}(v,t)\Big|_{v=1} = \left\{
\begin{array}{clr}
\ds \sum_{r=1}^t {\cal B}_r(\psi_1,\ldots,\psi_r;1), &  \ds t  = 1,2,\ldots \\ \\
  \ds 0, & \ds t = 0. &
 \end{array}\right.
\eeq
\section{Expected sojourn times on sites}
\label{recurrence_transience_ADTRW}
An important issue in random walk theory are recurrence/transience features. If the walker in an infinitely long walk returns 
to the departure site with probability one, the walk is said `recurrent' and
'transient' if this probability is smaller than one.
A vast literature on this topic exists \cite{SpitzerF1976,Feller1971,LawlerLimick2012}.
For the simple symmetric walk on $\mathbb{Z}^d$ the celebrated recurrence theorem
was established by P\'olya \cite{PolyaG1921}.
Recurrence and transience for symmetric L\'evy flights in multi-dimensional lattices and fractal features in  distributions were analyzed in \cite{TMM-APR-ISTE2019} and a recurrence theorem for these motions was established \cite{TMM-APR2017-recurrence,TMM_APR-recurrence2018}.  Further models considering recurrence and transience for modified L\'evy motions emerged only recently \cite{PagniniVitali2021}.

Recurrence/transience of a walk is an intrinsic property linked to the
{\it expected sojourn time (EST)} of an infinitely long walk.
The EST on site $n$ (departure site $m$) in an infinitely long ADTRW (i.e. $t\to \infty$) can be extracted from the generating function of the transition matrix (\ref{evolution}):
\beq
\label{tau_sojourn}
\begin{array}{clr}
\ds \mathbb{E}[\tau_{m,n}]  & =\mathbb{E}[ \tau_{0,n-m} ] = \ds \sum_{t=0}^{\infty} [{\mathbf P}(t)]_{m,n} =
[{\bar \Lambda}({\mathbf W}^{+},{\mathbf W}^{-},u)]_{m,n}\big|_{u=1}  & \\ \\
& =  \ds \frac{1}{\pi} \Re \int_0^{\pi}
\frac{W^{-}(\varphi)[{\mathbf 1}-{\bar \psi}(W^{-}(\varphi) )]}{(1-W^{-}(\varphi))\{W^{-}(\varphi)-W^{+}(\varphi){\bar \psi}[W^{-}(\varphi)]\}} e^{i\varphi(n-m)} {\rm d}\varphi. &
\end{array}
\eeq
Here $\Re$ extracts the real part of the following complex quantity.
\\[1mm]
The ADTRW is transient if $\ds \mathbb{E}[ \tau_{0,0} ] < \infty $ and recurrent if this quantity diverges. It is sufficient to consider the EST
on the departure site in order to verify recurrence/transience and we have that
$\ds \mathbb{E}[ \tau_{m,n}] / \mathbb{E}[ \tau_{0,0}] \leq 1 $.
In a transient walk a site is visited only a finite number of times as $t \to \infty$ whereas in a recurrent walk ($ \mathbb{E}[ \tau_{0,0} ] = \infty $) infinitely often by recurrent visits. 
\subsection{Recurrence/Transience features of the simple ADTRW}
\label{recurrence_transience}
The goal of the present part is to explore recurrence/transience features of simple ADTRWs
by means of analyzing their EST.
The EST on site $n$ (departure site $0$) of
Eq. (\ref{tau_sojourn}) has the canonical form
\beq
\label{gen_trans_feature}
 \mathbb{E} [ \tau_{0,n} ]_{simple} =\frac{1}{\pi} \Re \int_0^{\pi}e^{i n\varphi}{\bar \Lambda}(e^{-i\varphi},e^{i\varphi},u)\bigg|_{u=1}{\rm d}\varphi ,\hspace{0.5cm}
 {\bar \Lambda}(e^{-i\varphi},e^{i\varphi},u) = {\bar {\cal P}}(e^{-2i\varphi},u e^{i\varphi}), \hspace{0.3cm} n \in \mathbb{Z}.
\eeq
with transition matrices
 ${\mathbf W}_1^{+}={\hat T}_{-1}{\mathbf 1}$, ${\mathbf W}_1^{-}={\hat T}_{1}{\mathbf 1}$ having eigenvalues $W_1^{+}(\varphi)=e^{-i\varphi}$ and  $W_1^{-}(\varphi)=e^{i\varphi}$ ($\varphi \in (-\pi,\pi]$) due to the occurrence of (positive and negative) unit jumps. For what follows it is important to keep in mind that the transition matrices all have the same set of eigenvectors: The left eigenvectors are row-vectors with components $e^{i\varphi n}/\sqrt{2\pi}$ (to fulfill
$ {\hat T}_{-1}e^{i\varphi n}/\sqrt{2\pi}=\sum_{s=-\infty}^{\infty}\delta_{s,n-1}e^{i\varphi s}/\sqrt{2\pi} = e^{-i\varphi}e^{i\varphi n}/\sqrt{2\pi}$), and the right eigenvectors are column-vectors with components $e^{-i\varphi m}/\sqrt{2\pi}$, see Eq. (\ref{canocic}).
Note that ${\bar \Lambda}(1,1,u) =1/(1-u)$
(see Eq. (\ref{similar_gen})) is diverging for $u\to 1$. Correspondingly the generating function ${\bar \Lambda}(e^{-i\varphi},e^{i\varphi},1)$ becomes singular for $\varphi \to 0$ which is the only singularity in the integration interval. The type of this singularity is crucial 
for the integrability of Eq. (\ref{gen_trans_feature})
at $\varphi=0$ and hence to understand whether the walk is recurrent or transient. 
In order to explore this singular behavior we expand 
\beq
\label{small_varphi_exp}
\Re\left\{ {\bar {\cal P}}(e^{-2i\varphi},e^{i\varphi})\right\} = \Re\left\{
\frac{[1-{\bar \psi}(e^{i\varphi})]}{(1-e^{i\varphi})}\frac{1}{[1-e^{-2i\varphi}{\bar \psi}(e^{i\varphi})]}\right\}
\eeq
for small $\varphi$. To this end we need to account for the following possible cases: the interarrival time density $\psi(t)$ is either (a) `fat-tailed' (FT) or (b) `light-tailed' (LT). To capture this feature we expand the generating function ${\bar \psi}(u)={\bar \psi}_{\mu}(u)$ around the critical value $u=1$
which gives
\beq
\label{expansion}
{\bar \psi}_{\mu}(u) =  1 -  A_{\mu} (1-u)^{\mu} + o((1-u)^{\mu})  ,\hspace{0.5cm} \mu \in (0,1]
\hspace{0.5cm} |1-u| \to 0.
\eeq
containing the positive constant $A_{\mu} >0$ (independent of $u$) and $\mu$ indicates the lowest order occurring in this expansion.
Let us now analyze the FT and LT cases separately.
\\[3mm]
\noindent {\bf (a)} $\mu \in (0,1)$:\, ${\bar \psi}_{\mu}(t)$ fat-tailed (FT):
\\
We call a waiting-time density `fat-tailed' (FT) if expression (\ref{expansion}) is weakly singular at $u=1$ leading to an asymptotic power-law decay
 $\psi_{\mu}(t) \sim A_{\mu}(-1)^{t-1}\binom{\mu}{t} \sim \frac{A_{\mu}\mu}{\Gamma(1-\mu)}t^{-\mu-1} $ ($t \to \infty$) which is of the same type as in the Sibuya distribution. Therefore, the Sibuya distribution is a prototypical example for a FT distribution and of utmost importance. We will consider it closely in Section \ref{Sibuya_walk}.
 Eq. (\ref{small_varphi_exp}) then takes for $\varphi$ small 
\beq
\label{small_varphi_asymp}
\begin{array}{clr}
\ds {\bar {\cal P}}_{\mu} & \ds \sim
\frac{ A_{\mu}(1-e^{i\varphi})^{\mu-1}}
{1-e^{-2i\varphi}[1-A_{\mu}(1-e^{i\varphi})^{\mu}]} \sim \frac{i}{\varphi}\frac{1}{(1+\frac{2}{A_{\mu}}i^{\mu+1}\varphi^{1-\mu})},  \hspace{1cm} (\varphi \to 0) & \\ \\
 & \ds \sim \frac{i}{\varphi}+\frac{2 \varphi^{-\mu}}{A_{\mu}}i^{\mu}.  &
\end{array}
\eeq
Taking the real part shows the weakly singular behavior
\beq
\label{real_part}
\Re\{{\bar {\cal P}}_{\mu}\} \sim \frac{2\varphi^{-\mu} }{A_{\mu}}\cos{\left(\frac{\mu\pi}{2}\right)} >0,  \hspace{0.5cm} \mu \in (0,1) ,\hspace{0.5cm} (\varphi \to 0+)
\eeq
and hence it is integrable at $\varphi=0$.
We conclude that in the fat-tailed range $\mu \in (0,1)$ the
integral (\ref{gen_trans_feature}) exists, i.e. the EST on the sites is finite.
Therefore, simple
ADTRWs with generator processes of fat-tailed interarrival time densities always are transient.
A sufficient criteria is the weakly singular behavior
of the state probability generating functions
\beq
\label{fat-tail_survival}
\sum_{t=0}^{\infty}
{\mathbb P}(N(t)=n) =  \lim_{u\to 1}\frac{1-{\bar {\psi}(u)}}{1-u} {\bar \psi}(u)^n = \lim_{u\to 1} {\bar \Phi}^{(0)}_{\mu}(u) \sim A_{\mu}(1-u)^{\mu-1}
\to \infty,\hspace{0.3cm}
\, \mu \in (0,1)
\eeq
independent of state $n$
reflecting the universal asymptotic power-law behavior \\ $\Phi^{(n)}_{\mu}(t) \sim  A_{\mu}(-1)^t \binom{\mu-1}{t}\sim A_{\mu} \frac{t^{-\mu}}{\Gamma(1-\mu)}$ ($t \to \infty$), see e.g.\ \cite{MichelitschPolitoRiascos2021}.
\\[3mm]
\noindent {\bf (b)} $\mu=1$:\, ${\bar \psi}_{1}(t)$ light-tailed (LT):
\\
We call an interarrival time density `light-tailed' (LT) if its decay for large $t$ is at least geometrical
or faster, i.e.\ there are constants $C,\xi>0$ such that 
$|{\bar \psi}(t)| \leq C e^{-t \xi}$ for $t \to \infty$. As a consequence LT densities have finite moments (see also
Appendix \ref{Appendix_complex}).
For our convenience we introduce the complex variable $z=e^{i\varphi}$ to rewrite the EST on the sites $n \in \mathbb{Z}$ in expression (\ref{gen_trans_feature}) for an infinitely long walk as a closed complex contour integral over the unit circle $|z|=1$, namely 
\beq
\label{complex_int}
\begin{array}{clr}
\ds \mathbb{E}_{simple}[ \tau_{0,n} ] & = \ds \lim_{\epsilon \to 0+}{\bar \Lambda}({\hat T}_{-1},{\hat T}_{1},e^{-\epsilon})\delta_{0,n} = \lim_{\epsilon \to 0+}
\sum_{t=0}^{\infty} e^{-t\epsilon} [P_{0,n}(t)]_{simple} , \hspace{1cm}\, (\epsilon >0) & 
\\ \\
& = 
\ds \lim_{\epsilon \to 0+} \ds  \frac{1}{2\pi}\int_{-\pi}^{\pi} e^{i n\varphi}\,
{\bar {\cal P}}(e^{-2i\varphi},e^{i(\varphi+i\epsilon)})\,{\rm d}\varphi   &  \\ \\ &=\ds \lim_{\epsilon \to 0+} \ds
\frac{1}{2\pi i}\oint_{|z|=1}\frac{[1-{\bar \psi}( z e^{-\epsilon})]}{(1-z e^{-\epsilon})}
\frac{z^{n+1}}{[z^2-{\bar \psi}(ze^{-\epsilon})]}{\rm d}z & 
\\ \\
  &=: \ds P.V.\,
\frac{1}{2\pi i}\oint_{|z|=1}
\frac{z^{n+1} \, {\bar \Phi}^{(0)}(z)}{(z-1)[z+1 -{\bar \Phi}^{0}(z)]}{\rm d}z &
\\ \\ & = \ds P.V.\,
\frac{1}{2\pi i}\oint_{|z|=1}
\frac{z^{n} \, {\bar \Phi}^{(0)}(z)}{z-{\bar g}(z)}{\rm d}z.          & 
\end{array}
\eeq
In the last line we introduced ${\bar \psi}(z)=z{\bar g}(z)$ with the auxiliary generating function 
${\bar g}(z)$ (see Appendix \ref{Appendix_complex} for essential features). 
The integrand of Eq. (\ref{complex_int}) is singular at $z_1=1$ and  corresponds to the singularity of the quantity (\ref{small_varphi_exp}) at $\varphi=0$. In order to achieve a regularization we consider instead of $u=1$ the limit $u=e^{-\epsilon}$ for infinitesimally positive $\epsilon$ with
${\bar \Lambda}(e^{-i\varphi},e^{i\varphi},u)\big|_{u=e^{-\epsilon}} = {\bar {\cal P}}(e^{-2i\varphi},e^{i(\varphi+i\epsilon)})$. In this way we shift the singularity at $z_1=1$ infinitesimally away from the unit circle and obtain a well defined integral.
\\[1mm]
Let us explore in which direction this regularization procedure shifts the singularity $z_1=1$.
First, we observe as a consequence of the LT feature that the generating function of the survival probability ${\bar \Phi}^{(0)}(z)$ is analytical and finite at $z=1$ (see Eqs. (\ref{gener_function_expansion}), (\ref{survival_property})) thus the singularity at $z_1=1$ is due to the zero $z^2-{\bar \psi}(z)=0$.
To this end we put the shifted zero to $z_1(\epsilon)=e^{a\epsilon}$ which is infinitesimally close to one, where the constant $a$ is independent of $\epsilon$ and has to be determined 
from $z^2-{\bar \psi}(ze^{-\epsilon})=0$, leading to
\beq
\label{D1_epsilon}
e^{2a\epsilon}-{\bar \psi}(e^{(a-1)\epsilon}) = 0.
\eeq
The first order in $\epsilon$ must identically vanish, which yields
\beq
\label{leading_to_a}
a=\frac{A_1}{A_1-2} = -\frac{1+g_1}{1-g_1} ,\hspace{1cm} g_1=\frac{d}{dz}{\bar g}(z)|_{z=1}=A_1-1
\eeq
where the constant $A_1=\frac{d}{dz}{\bar \psi}(z)|_{z=1}= \sum_{t=1}^{\infty}t\psi(t)$ (with $A_1 \geq 1$) indicates the expected waiting time between successes (positive unit jumps).
The sign of the parameter $a$ is crucial in order to see whether $z_1=e^{a\epsilon}$ is within or outside the unit disc.
It follows that (i) $a>0$ for $A_1>2$ and (ii) $a<0$ for $A_1<2$, and therefore
\beq
\label{z_one_shift}
\begin{array}{clr}
z_1 & \sim  e^{\frac{\epsilon A_1}{A_1-2}} \sim 1 + , \hspace{1cm}{\rm for}\hspace{1cm} A_1> 2, & \\[2mm]
 z_1 & \sim  e^{\frac{\epsilon A_1}{A_1-2}} \sim 1- , \hspace{1cm}{\rm for}\hspace{1cm} 1\leq A_1 < 2. &
 \end{array}
\eeq
Hence the residue at $z_1=1$  contributes to contour integral
(\ref{complex_int}) for $1 \leq A_1<2$ but does not contribute in the range $A_1>2$. 
Interestingly, the sign of the infinitesimal shift is solely determined by the sign of $B=2-A_1$. Later on we will see more closely that $B$ is a measure for `bias' in an asymptotic sense which emerges in a simple ADTRW for large $t$.
\\[1mm]
Crucial for the further analysis is the expansion of the generating function of the survival probability
\beq
\label{gener_function_expansion}
\begin{array}{clr}
\ds {\bar \Phi}^{(0)}(z) &= \ds \frac{1}{z-1}\left(-1 + {\bar \psi}(1) + \sum_{\ell=1}^{\infty}\frac{(z-1)^{\ell}}{\ell !}\frac{d^{\ell}}{dz^{\ell}}{\bar \psi}(z)\big|_{z=1}\right) & \\ \\ & = \ds A_1+\sum_{\ell=2}^{\infty}A_{\ell}(z-1)^{\ell-1}, & \hspace{-3cm}\ds A_{\ell} = \frac{1}{\ell!} \frac{d^{\ell}}{dz^{\ell}}{\bar \psi}(z)\big|_{z=1} \geq 0
\end{array}
\eeq
which is analytic on the unit disc $|z|\leq 1$ where all $A_{\ell}$ are finite as a consequence of the LT feature of $\psi(t)$. The lower bound of the expected waiting time is $A_1=1$ and occurs only in the trivial case when each trial almost surely is a success, corresponding to the interarrival time density
$\psi_{trivial}(t)=\delta_{1,t}$ with ${\bar \psi}_{trivial}(z)=z$.
Further, we observe that
${\bar \Phi}^{(0)}(z)\big|_{z=0}=A_1 -\sum_{l=2}^{\infty} 
(-1)^{\ell} A_{\ell}=1$ recovering the initial condition of the survival probability.
On the other hand $\frac{d}{dz}{\bar g}(z)\big|_{z=1}=g_1=A_1-1 \geq 0$ thus $g_1=\sum_{l=2}^{\infty} 
(-1)^{\ell} A_{\ell}$.
In particular, we have that
\beq
\label{survival_property}
{\bar \Phi}^{(0)}(z)\big|_{z=1}= {\bar \Phi}^{(n)}(z)\big|_{z=1} =\lim_{z\to 1}\frac{1-{\bar \psi}(z)}{1-z} =\frac{d}{dz}{\bar \psi}(z)\big|_{z=1} = \sum_{t=1}^{\infty} t \psi(t) = A_1 \geq 1
\eeq
yielding the expected interarrival time between successive successes.
We hence observe the inequality $1 \leq {\bar \Phi}^{(0)}(z)\leq A_1$ for $z \in [0,1]$ where we also use
that ${\bar \Phi}^{(0)}(z)$ is absolutely monotonic (AM) in that interval (see Appendix \ref{Appendix_complex}).
Clearly by using the LT feature (\ref{survival_property}) it follows for the state probabilities 
\beq
\label{convergence_gen}
\sum_{t=0}^{\infty}{\mathbb P}(N(t)=n) = \left({\bar \Phi}^{(0)}(z) ({\bar \psi}(z))^n\right)\big|_{z=1}=A_1 , \hspace{1cm} \forall n \in \mathbb{N}_0.
\eeq
Hence, recurrence/transience solely depends on the singularities of the part $\frac{z}{(z^2-{\bar \psi}(e^{-\epsilon} z))}$ ($\epsilon \to 0+$)  within the unit disc where
\beq
\label{zeros}
D(z) =  z^2-{\bar \psi}(z) = - z \left({\bar g}(z)-z\right)  ,\hspace{3cm} {\bar \psi}(z)=z {\bar g}(z).
\eeq
We prove in Appendix \ref{Appendix_complex} that the complex function
 ${\bar g}(z)-z$ has a canonical representation of the form

\beq
\label{we_have_A1sm2}
{\bar g}(z)-z = (z-1)(z-r) e^{h(z)} , \hspace{0.5cm} r \in \mathbb{R}^+, \hspace{0.5cm} |z| \leq 1,
\eeq
which exists {\it at least} on the unit disc $|z|\leq 1$. 
The zero $r$ is real with the properties $r=r(A_1)>1$ for $A_1<2$ and $r(A_1)<1$ for $A_1>2$. Further, $r=1$ for $A_1=2$
which is the recurrent limit where the multiplicity of the zero $z=1$ then is two. The contribution $e^{h(z)}$ has no zeros and $h(z)$ is analytic at least on the unit disc $|z|\leq 1$.
Consult Appendix \ref{Appendix_complex} for a detailed discussion with proofs and the example of Poisson distributed waiting times where a canonical form (\ref{we_have_A1sm2}) for all $z\in \mathbb{C}$ exists (see Eq. (\ref{canonic_function-rep2})).
Due to the structure of the zeros of Eq. (\ref{we_have_A1sm2}) in the evaluation of the contour integral (\ref{complex_int}) by residue theorem we have to distinguish the two cases $1\leq A_1<2$ and $A_1>2$.
\\[1mm]
For $1\leq A_1=g_1+1<2$ ($g_1=\frac{d}{dz}{\bar g}(z)|_{z=1} \in (0,1)$) we get
 \beq
 \label{residuum}
 \mathbb{E}[ \tau_{0,n}]_{simple}  = \left\{\begin{array}{l} \ds  \frac{{\bar \Phi}^{0}(z)}{1-\frac{d}{dz}{\bar g}(z)}\big|_{z=1-} =\frac{1+g_1}{1-g_1}  = \frac{A_1}{2-A_1},  
  \hspace{2.5cm}
   n \in \mathbb{N}_0 \\ \\ \ds 
  \frac{A_1}{2-A_1}
    + \frac{1}{(|n|-1)!}\frac{d^{|n|-1}}{dz^{|n|-1}}\left(\frac{{\bar \Phi}^{(0)}(z)}{z-{\bar g}(z)}\right)\big|_{z=0}
    , \hspace{1cm} -n \in \mathbb{N}
   \end{array}\right.
 \eeq
 where we accounted for the properties (\ref{z_one_shift}), i.e.\ the fact that singularity at $z=1 \to 1-$ is for $A_1<2$ infinitesimally shifted inside the unit disc and therefore contributes whereas the second zero $r(A_1)>1$ does not contribute.
 The second line in Eq.\,(\ref{residuum}) refers to the sites $n<0$ where due to $z^{-|n|}$ an additional singularity at $z=0$ occurs.
 Remarkably, all sites $n \geq 0$ are 
 equally long visited in an infinitely long walk, namely $\mathbb{E}[ \tau_{0,n}]_{simple}=\mathbb{E}[ \tau_{0,0}]_{simple}= 
 \frac{1+g_1}{1-g_1} \geq 1$ with $0 \leq g_1 <1$ in the range $1 \leq A_1 <2$.
 Moreover, in the trivial strictly increasing walk where each jump almost surely is of size $+1$ with
 $A_1=1$ we have the minimum value $\tau_{0,n}=1$ $\forall n >0$ where the EST on all these sites is one.
We observe that $\mathbb{E}[ \tau_{0,n} ] \to \infty$ for $A_1 \to 2$ which is hence the recurrent limit.
In the recurrent limit the expected waiting time between
jumps of the same direction is $A_1=2$ which physically means that for $t\to \infty$ an equal expected numbers of jumps $+1$ and $-1$ occur. We expect therefore that a recurrent simple ADTRW is unbiased at least in an asymptotic sense, i.e. for an infinitely long observation time. We prove this assertion subsequently in this section (see asymptotic relation (\ref{simple_asymptotics}))
and explore the interplay of bias and recurrence/transience features.
To this end we consider the quantity $B=2-A_1$ which will turn out to contain crucial information on the bias.
\\[1mm]
 We call a simple ADTRW {\it strictly unbiased} if the following conditions (i) and (ii) are fulfilled. 
 \\[1mm]
 \noindent ({\bf i}) $A_1=2$ ($B=0$), i.e.\ the walk is recurrent and the expected position (\ref{simple_asymptotics}) is null in the limit of an infinitely long observation $t\to \infty$. 
\\[2mm]
\noindent ({\bf ii}) The state polynomial (\ref{similar}) and its generating function fulfill the exchange symmetry property ${\bar \Lambda}(a,b,u)={\bar \Lambda}(b,a,u)$, and as a consequence the transition matrix then is symmetric (self-adjoint) with
$$\Lambda({\hat T}_{-1},{\hat T}_1,t) = [\Lambda({\hat T}_{-1},{\hat T}_1,t)]^{\dagger} =\Lambda({\hat T}_{1},{\hat T}_{-1},t) $$
as the unitary shift operators $({\hat T}_{-1})^{\dagger}={\hat T}_1$ are adjoint to each other.
As a consequence of (ii)
the eigenvalues of the transition matrix 
$\Lambda(e^{-i\varphi},e^{i\varphi},t)= \Lambda(e^{i\varphi},e^{-i\varphi},t)$ (and of its generating function matrix ${\bar \Lambda}(e^{-i\varphi},e^{i\varphi},u)$)
 are real and even functions of $\varphi$. We also see that in this case the expected position of the walker (\ref{simple_walk_bias}) 
 \beq
 \label{expected_position_is_null}
 \mathbb{E}[Y_t]_{simple} = \left(\frac{\partial }{\partial a}- \frac{\partial }{\partial b} \right) \Lambda(a,b,t)\Big|_{a=b=1} = 0 ,\hspace{1cm} \forall t \in \mathbb{N}_0
 \eeq
 is null for all times as a consequence of the exchange symmetry
 $\Lambda(a,b,t)= \Lambda(b,a,t)$
 (where $0$ is the departure site). The occurrence of a symmetry in $\Lambda(a,b,t)$ such as the exchange symmetry $a \leftrightarrow b $ reflects a conserved quantity, namely $\mathbb{E}[Y_t]_{simple}=0$. We point out that this observation has a remarkable analogy with Noether's theorem which roughly tells us that each symmetry corresponds to a conserved quantity \cite{Noether1918}.
 If (ii) is fulfilled (i) is fulfilled, conversely if (i) is fulfilled (ii) does not necessarily hold true.
 To see this, consider the generating function (see Eq.\,(\ref{gener-limit}) below) of the expected position of the walker.
For a strictly unbiased simple ADTRW this quantity must vanish which yields a condition for the generating
function ${\bar \psi}_{unbiased}(u)$ of the generator process for which the walk is strictly unbiased, namely\footnote{$u/(1-u)^2$ being the generating function of $t$ and $\frac{\partial}{\partial v}{\bar P}(v,u)\big|_{v=1} =\frac{{\bar \psi}(u)}{(1-u)(1-{\bar \psi}(u))}$.}
\beq
\label{condition_strictly_unbiased}
 \frac{2{\bar \psi}_{unbiased}(u)}{(1-u)(1-{\bar \psi}_{unbiased}(u))} -\frac{u}{(1-u)^2} = 0 .
\eeq
This yields
\beq
\label{unbiased_generating-function}
{\bar \psi}_{unbiased}(u) = \frac{u}{2(1-\frac{u}{2})}
\eeq
that is the generating function ${\bar \psi}_B(u)=\frac{pu}{1-qu}$ with geometrically distributed waiting times $\psi_{B}(t)=pq^{t-1}$ of the {\it Bernoulli generator
process} for the symmetric case $p=q=1/2$. In particular, we have the state polynomial $\Lambda_B(a,b,t)=(a+b)^t/2^t$ (see Eq.\,(\ref{Lambda_bern})) which indeed fulfills the claimed exchange symmetry
$\Lambda_B(a,b,t)=\Lambda_B(b,a,t)$, i.e.\ condition (ii). 
We have for this recurrent walk indeed $A_1=1/p=2$, i.e.\ condition (i) also holds true. Since the result (\ref{unbiased_generating-function}) is unique for simple walks, it follows that the simple ADTRW with symmetric Bernoulli generator process with $p=q=1/2$ indeed is the only one which is {\it strictly unbiased}, i.e. fulfills conditions (i) and (ii) with Eq. (\ref{expected_position_is_null}) for all times $t\in \mathbb{N}_0$.
Further, note that this is consistent with the fact that the interarrival times between two consecutive successes and those between to consecutive fails are equally distributed.
 We consider briefly the Bernoulli ADTRW at the end of this section and consult Appendix 
 \ref{Appendix_C} as well as the references \cite{SpitzerF1976,Feller1971,PolyaG1921,RednerS,LawlerLimick2012}.
The class of simple walks with $A_1=2$ which do not have Bernoulli generator process are recurrent fulfilling (i) but they are not strictly unbiased since they do not fulfill (ii). The class of these simple walks is unbiased in an asymptotic sense, i.e. $\lim_{t\to \infty} \mathbb{E}[Y_t]_{simple} =0$ (see (\ref{simple_asymptotics}) for $A_1=2$). 
\\[1mm]
We devote the subsection \ref{simple_prescribed} to consider the class of admissible prescribed functions for the expected position $\mathbb{E}[Y_t]_{simple} = f(t)$ in an ADTRW.
\\[3mm]
The class of simple ADTRWs with $A_1\neq 2$ are both {\it biased} and transient (with finite EST $\mathbb{E}[ \tau_{0,n}]_{simple} \leq  \mathbb{E}[ \tau_{0,0}]_{simple} <\infty $). The walk is `right-biased' for $A_1<2$ ($B>0$) where the expected number of positive jumps dominates, and `left-biased' for $A_1>2$ ($B<0$) with domination of the expected number of negative jumps in infinitely long walks. 
 In the picture of Gambler's Ruin Problem the simple ADTRW defines in the range $B=2-A_1>0$ ($B\in (0,1]$) a long-time `winning strategy'. We will come back to the asymptotic behavior subsequently.
 \\[3mm]
 For $A_1>2$ ($g_1=\frac{d}{dz}{\bar g}(z)|_{z=1} \in (1,\infty)$) the EST integral (\ref{complex_int}) yields
 \beq
 \label{A_1_greater2}
 \mathbb{E}[ \tau_{0,n}]_{simple} =
 \left\{\begin{array}{l} \ds  
 \frac{r^{n}\, {\bar \Phi}^{(0)}(r)}{1-\frac{d}{dz}{\bar g}(z)\big|_{z=r}}  \hspace{1cm}
   n \in \mathbb{N}_0 \\ \\ \ds 
     \frac{r^{n}\, {\bar \Phi}^{(0)}(r)}{1-\frac{d}{dz}{\bar g}(z)\big|_{z=r}}
    + \frac{1}{(|n|-1)!}\frac{d^{|n|-1}}{dz^{|n|-1}}\left(\frac{{\bar \Phi}^{(0)}(z)}{z-{\bar g}(z)}\right)\big|_{z=0}
    , \hspace{1cm} -n \in \mathbb{N}.
    \end{array}\right.
 \eeq
 where the zero of (\ref{we_have_A1sm2}) $r=r(A_1) \in (0,1)$ and the singularity at $z=1 \to 1+$ is here infinitesimally shifted outside the unit disc and therefore does not contribute (see Eq.\,(\ref{z_one_shift})). By using the convex feature 
 of ${\bar g}(z)$ we show in Appendix \ref{Appendix_complex} that $\frac{d}{dz}{\bar g}(z)\big|_{z=r} \in (0,1)$ (see Eq.\,(\ref{derivative_in_r})) ensuring that the quantity (\ref{A_1_greater2}) is strictly non-negative.
 The EST on the sites $n\geq 0$ on the right of the departure node decays geometrically
 as $r^n$ which is consistent with the picture that this walk with  $B=2-A_1 <0$ is left-biased for $t\to \infty$. Therefore 
 the geometric decrease of the EST for positive sites physically makes sense.
\\[1mm]
Indeed, expressions
(\ref{residuum}) and (\ref{A_1_greater2})
are both positively singular for the recurrent limit
 $A_1=2$ ($B=0$), thus integrability of the contour integral (\ref{complex_int}) breaks down as the multiplicity of zero 
 $z=r=1$ in (\ref{we_have_A1sm2}) is two.
 The divergence in formula (\ref{A_1_greater2}) for
 $A_1=1+g_1 \to 2+$ can be seen by $r(A_1) \to 1-$ and accounting for $\frac{d}{dz}{\bar g}(z)\big|_{r(A_1)} \to 1-$ in this limiting case (see Appendix \ref{Appendix_complex} for a detailed outline of related properties).
\\[1mm]
Consider again the EST (\ref{residuum}) in the range $A_1<2$ for the left neighbor site $n=-1$ 
from the departure site. We have
\beq
\label{sojournm1}
\mathbb{E}[\tau_{0,-1}]_{simple} =  \frac{A_1}{2-A_1} -\frac{{\bar \Phi}^{0}(z)}{{\bar g}(z)}\big|_{z=0} = 
\frac{A_1}{2-A_1} - \frac{1}{\alpha_1} , \hspace{0.5cm} 1\leq A_1<2
\eeq
where we used ${\bar g}(0)=\psi(1)=\alpha_1$ together with the initial condition of the survival probability ${\bar \Phi}^{(0)}(z)\big|_{z=0}=\Phi^{(0)}(t)\big|_{t=0}=1$. We see that 
$\mathbb{E}[\tau_{0,-1}]_{simple} < \mathbb{E}[\tau_{0,0}]_{simple}$ which clearly reflects the fact that the walk with $B=2-A_1>0$ is right-biased as $t\to \infty$ where the jumps in positive direction dominate. For a proof of non-negativeness
of $\mathbb{E}[ \tau_{0,-1}]_{simple}$ consult Appendix \ref{Appendix_complex}.
For the trivial walk with $A_1=1$, $\alpha_1=1$ we have $\mathbb{E}[ \tau_{0,0}]=1$, i.e.\ the walker is present on the departure site almost surely only during one time unit following its departure at $t=0$.
On the other hand we then have
 $\mathbb{E}[ \tau_{0,-1}] =0$, i.e.\ the site $-1$ is almost surely not visited in the strictly increasing trivial walk.
 Consider the trivial walk $A_1=1$ for all $n<0$ which yields (where in this case ${\bar \Phi}^{0}(z)=1$ and ${\bar g}(z)=1$)
 \beq
 \label{trivial_case}
 \mathbb{E}[\tau_{0,-n}]_{simple} =1 +\frac{1}{(n-1)!}\frac{d^{n-1}}{dz^{n-1}}\frac{1}{(z-1)}\big|_{z=0}=1+\frac{(-1)^{n-1}}{(z-1)^n}
 \big|_{z=0} = 0 \hspace{1cm} n \in \mathbb{N}
 \eeq
i.e.\ all sites on the negative side of the departure node are almost surely not visited. This is perfectly consistent with the physical picture of the trivial strictly increasing walk performing unit jumps $+1$ almost surely in each time increment.
\\[1mm]
A further interesting quantity also is the probability that a site $n$ in an infinitely long walk is ever visited (for $n=0$ that the walker ever returns to the departure site).
This quantity is related with the EST by \cite{Feller1971,PolyaG1921,TMM-APR2017-recurrence,TMM_APR-recurrence2018}
\beq
\label{ever_prob}
F_{0,n} = \frac{\mathbb{E}[\tau_{0,n}]_{simple} -\delta_{0,n}}{\mathbb{E}[\tau_{0,0}]_{simple}}
\eeq
and yields for $1 \leq A_1 \leq 2$: 
\beq
\label{trnasience}
\begin{array}{clr}
\ds F_{0,0} & = \ds  \frac{2(A_1-1)}{A_1} & \\ \\
\ds F_{0,n} & = \ds  1 , & \ds  n=\{1,2,\ldots\} \in \mathbb{N}. 
\end{array}
\eeq
The quantity $1-F_{0,0}=(\mathbb{E}[\tau_{0,0}]_{simple})^{-1}$ can be 
interpreted as the `escape probability', i.e. 
the probability that the walker never returns to the departure site. 
Further, we have
$0\leq F_{0,0} \leq 1$ where in the recurrent limit $F_{0,0}=1$ 
and we have $F_{0,0} <1$ in all transient (biased) cases ($A_1\neq 2$).
 For $A_1=1$ have $F_{0,0}=0$ (almost surely no return to the departure site) and a.s. no visits on negative sites $F_{0,-n}=0$ ($-n < 0$) in the trivial strictly increasing walk. In the recurrent limit $A_1=2$ we have for all sites $F_{0,0}=F_{0,n}=1$ ($n \in \mathbb{Z}$), i.e.\ each site is almost surely ever visited.
 \\[1mm]
 It appears instructive to consider here also the connection of bias and expected position of the walker (\ref{simple_walk_bias}) for large times $t\to \infty$. The generating function ${\bar Y}_{simple}(u)$ of this quantity is
 \beq
 \label{gener-limit}
 \begin{array}{clr}
 \ds {\bar Y}_{simple}(u) &= \ds 
 2\frac{\partial {\bar {\cal P}}}{\partial v}(v,u)|_{v=1} -\frac{u}{(1-u)^2} & \\ \\ 
  & = \ds  \frac{2{\bar \psi}(u)}{(1-u)(1-{\bar \psi}(u))} -\frac{u}{(1-u)^2} & 
 \end{array}
 \eeq 
 where $\frac{u}{(1-u)^2}= \sum_{t=1}^{\infty}tu^t$ is the generating function of $t \in \mathbb{N}_0$ and 
 ${\bar N}(u)=\frac{\partial}{\partial v}{\bar {\cal P}}(v,u)\big|_{v=1} $ of the expected number of arrivals. In order to capture the large time asymptotic behavior we expand 
 (\ref{gener-limit}), $u\to 1-$, and arrive at
 \beq
 \label{gen-result}
 {\bar Y}_{simple}(u)\sim \frac{(2-A_1)}{A_1} \frac{u}{(1-u)^2} ,\hspace{2cm} u \to 1-
 \eeq
 which gives the asymptotics for the expected position of the walker for $t$ large
 \beq
 \label{simple_asymptotics}
 \mathbb{E}[Y_t]_{simple} \sim \frac{(2-A_1)}{A_1} t ,\hspace{3cm} t \to \infty.
 \eeq
 The sign of this quantity defines the bias in an asymptotic sense where this formula holds for all $B$ and indeed in the recurrent case
 $B=0$ the walk is unbiased in the limit $t \to \infty$.
 For the Bernoulli generator process relation (\ref{simple_asymptotics}) recovers the well known classical result (\ref{mean}) (see Appendix \ref{Appendix_C}). For the trivial strictly increasing walk with $A_1=1$ we have necessarily
 $ \mathbb{E}[Y_t]_{simple} =t $. On the other hand the asymptotic behavior $\mathbb{E}[Y_t]_{simple} = -t$ is approached in the fat-tailed
 limit $A_1 \to \infty$ and is indeed the dominant contribution in the asymptotic formula (\ref{as_large_t}) for the simple Sibuya ADTRW considered in Section \ref{Sibuya_walk}.
 \\[3mm]
 \noindent {\it Simple Bernoulli ADTRW}
 \\
 Let us compare some of these results with the 
 case of the Bernoulli generator process (see also Appendix \ref{Appendix_C}).
 The Bernoulli trial process has geometric light-tailed waiting-time density
 $\psi_{B}(t)=pq^{t-1}$ ($p+q=1$) with generating function
 ${\bar \psi}_B(z)=\frac{pz}{1-qz}$ where
 $A_1= \frac{d}{dz}\frac{pz}{1-qz}|_{z=1}= 1/p $ and we have
 for the canonical form (\ref{we_have_A1sm2})
\beq
\label{poles}
{\bar g}_{B}(z)-z= \frac{p}{1-qz} -z = (z-1)(z-\frac{p}{q})\, \frac{q}{1-qz} =(z-1)(z-\frac{p}{q})e^{h_B(z)}
\eeq
with the zeros $1,r_B=p/q$ and $h_B(z)=\log(q)-\log(1-qz)$.
For $A_1=1/p<2$, i.e. $p>q$ the second zero $r_B>1$ is outside the unit disc in agreement with our above result that for $A_1<2$ the function ${\bar g}_B(z)-z=0$ has the only (infinitesimally shifted) zero 
$z_1=1-$ in the unit disc (See Appendix \ref{Appendix_C} for details).
The second zero $r_B= p/q $
is within the unit disc only for $A_1>2$ ($p<q$), and outside for $1\leq A_1 <2$ ($p>q$) in agreement with the general behavior outlined above (and see Appendix \ref{Appendix_C}).
On the other hand $p=q=\frac{1}{2}$ ($A_1=2$) represents the recurrent limit and represents the only existing strictly unbiased simple ADTRW where the multiplicity of the zero
$z_1=r_B=1$ in expression (\ref{poles}) then is two and $\mathbb{E}[ (Y_t)_B]_{simple}=0$ $\forall t$.
\\[1mm]
One obtains with relations (\ref{residuum}) and (\ref{A_1_greater2})
for the EST on the departure site in an infinitely long simple Bernoulli ADTRW 
\beq
\label{bernoulli_sojounr}
\ds \mathbb{E}[ \tau_{0,0} ]_{B} = \ds  \lim_{\epsilon \to 0+} \frac{1}{2\pi i}\oint_{|z=1|}\frac{{\rm d}z}{(-q)e^{-\epsilon}[z^2-e^{\epsilon}\frac{z}{q}  +\frac{p}{q}]} = \frac{1}{|p-q|}
\eeq
which is a classical result given by Feller \cite{Feller1971} (see Chapter VIII).
\subsection{Prescribed expected position and bias in a simple ADTRW}
\label{simple_prescribed}
In many applications it may be interesting to prescribe in a simple ADTRW not the generator process, but the expected position of the walker $\mathbb{E}[Y_t]_{simple} =f(t)$ where $f(t)$ has to be an admissible function which fulfils $f(t)\big|_{t=0}=0$, $t \in \mathbb{N}$. 
 It follows from formula (\ref{simple_walk_bias}) that the class of admissible functions $f(t)$ is restricted by 
 $-t \leq f(t) \leq t$ and $f(t)=-t+2C(t)$ where $C(t)$ is non-negative and non-decreasing with $0\leq C(t) \leq t$ as a consequence of $0 \leq N(t) \leq t$ with
 $C(t) = \mathbb{E}[N(t)]$ and therefore is defined as well on $t\in \mathbb{N}$ (with $C(0)=0$). Its generating function ${\bar C}(u)$ is absolutely monotonic.
 Considering Eq. (\ref{simple_walk_bias}) we have that
 \beq
 \label{leads_to}
\mathbb{E}[Y_t]_{simple} = 2 \frac{\partial}{\partial v}{\cal P}(v,t)\Big|_{v=1} - t  = f(t) = 2C(t) -t ,\hspace{1cm} t \in \mathbb{N}_0.
  \eeq
  Let ${\bar f}(u)=\sum_{t=1}^{\infty}f(t)u^t$ be the generating function which we assume to converge at least for  $|u|<1$. It is convenient to put ${\bar f}(u)= u{\bar k}(u)/(1-u)^2$ (i.e. $f(t) = t \star {\bar k}(t)$) 
  and ${\bar k}(u)=2{\bar c}(u)-1$ with ${\bar C}(u)=u{\bar c}(u)/(1-u)^2$. 
  Taking then generating function on both sides of Eq. (\ref{leads_to}) it yields 
 \beq
 \label{psi_f}
 {\bar \psi}_f(u) = \frac{u(1+{\bar k}(u))}
 {2[1-\frac{u}{2}(1-{\bar k}(u)]} =  \frac{u{\bar c}(u)}{1-u +u{\bar c}(u)}   , \hspace{1cm} {\bar k}(u)=2{\bar c}(u)-1
 \eeq
where ${\bar c}(u) \in [0,1]$ for $u\in [0,1)$.
We observe that ${\bar \psi}_f(u)|_{u=0}=0$ and ${\bar \psi}_f(u)|_{u=1}=1$, i.e. $\psi_f(t)$ has the good properties of a waiting-time density supported on $\mathbb{N}$. 
Generally, the resulting generator process corresponding to Eq. (\ref{psi_f}) allows both LT and FT waiting-time densities $\psi_f(t)$.
\\[1mm]
Consider the linear law $f(t)=b_0t$ with  constant $b_0 \in [-1,1]$,  (i.e.  $k_{b_0}(t)=b_0\delta_{t,0}$).
Then Eq. (\ref{psi_f}) takes the form 
\beq
 \label{psi_f_b0}
 {\bar \psi}_{b_0}(u) = \frac{u(1+b_0)}
 {2[1-\frac{u}{2}(1-b_0]} 
\eeq
where we identify this expression with the generating function of the LT waiting time-density of the {\it Bernoulli generator process} with $p=(1+b_0)/2$ (and ${\bar c}(u)=p$ constant), $q=(1-b_0)/2$, i.e.\ $b_0=p-q$, recovering the well known relation 
(\ref{mean}) for the expected position of the walker. The case $f(t) = 0 \,\forall \, t \in \mathbb{N}_0$ is covered by
${\bar k}(u)=b_0=0$ in Eq. (\ref{psi_f_b0}) and recovers the generating function (\ref{unbiased_generating-function}) corresponding to the {\it strictly unbiased simple ADTRW}, i.e.\ the ADTRW with the symmetric Bernoulli generator process.
\section{Sibuya ADTRW}
\label{Sibuya_walk}
In this section we consider an ADTRW with generator process of Sibuya distributed interarrival times as a FT prototypical case. We refer this walk to as `Sibuya ADTRW'. 
The probabilities of first success (\ref{Sibuyaalpha}) in a Sibuya trial process have generating function

\beq
\label{Sibuya_gen}
{\bar \psi}_{Sibuya}(u)= \sum_{t=1}^{\infty}(-1)^{t-1}\binom{\beta}{t} u^t   = 1-(1-u)^{\beta} ,\hspace{0.5cm} |u|\leq 1, \hspace{0.5cm} \beta \in (0,1).
\eeq
The FT feature of the Sibuya distribution is reflected by the divergence of the expected interarrival time between successes: $\frac{d}{du}{\bar \psi}_{Sibuya}(u)\big|_{u=1} =\beta(1-u)^{\beta-1}|_{u=1} \to \infty$ (see also the asymptotic expansion
(\ref{expansion})).
The generating function (\ref{gen_state_pol}) of the Sibuya state polynomial yields
\beq
\label{Sibuya_state_pol_gen}
{\bar {\cal  P}}_{Sibuya}(v,u)= \frac{(1-u)^{\beta-1}}{1-v +v(1-u)^{\beta}} ,\hspace{0.5cm} |u| < 1 , \hspace{0.5cm} |v| \leq 1
\eeq
where necessarily
${\bar {\cal  P}}_{Sibuya}(1,u)=\frac {1}{1-u}$
(corresponding to the normalization of the Sibuya state probabilities) holds true.
Then we have for the generating function (\ref{similar_gen})
\beq
\label{related_gen}
   {\bar \Lambda}_{Sibuya}(a,b,u)   = {\bar {\cal  P}}_{Sibuya}\left(\frac{a}{b}, ub\right) =\frac{b(1-bu)^{\beta-1}}{b-a +a(1-bu)^{\beta}}
\eeq
with the limiting cases ${\bar \Lambda}_{Sibuya}(1,0,u) =\frac{1}{1-\beta u}$ and
${\bar \Lambda}_{Sibuya}(0,1,u)  =(1-u)^{\beta-1}$
(see Appendix \ref{Appendix_A} for details).
The Sibuya state polynomial is obtained from
\beq
\label{Sibuya_state}
{\cal P}_{Sibuya}(v,t) = \frac{1}{t!} \frac{d^t}{du^t} {\bar {\cal  P}}_{Sibuya}(v,u)\big|_{u=0} = \frac{1}{t!} \frac{d^t}{du^t} \left\{
\sum_{n=0}^{\infty} v^n (1-u)^{\beta -1}[1 - (1-u)^{\beta}]^n\right\}\bigg|_{u=0}
\eeq
and yields
\beq
\label{Sibuya-state_eval}
\begin{array}{clr}
\ds {\cal P}_{Sibuya}(v,t) & = \ds  \sum_{n=0}^{t} v^n{\mathbb P}(N_{Sibuya}(t)=n)= \frac{1}{t!} \frac{d^t}{du^t} \left\{
\sum_{n=0}^{t} v^n \sum_{\ell=0}^n \binom{n}{\ell}(-1)^{\ell} (1-u)^{\beta \ell +\beta -1}\right\}\big|_{u=0} & \\ \\
 &  =   \ds \frac{(-1)^t}{t!}\sum_{n=0}^{t} v^n \sum_{\ell=0}^n (-1)^{\ell} \binom{n}{\ell}  
 \frac{\Gamma(\beta[\ell+1])}{\Gamma(\beta[\ell+1] -t)}&  
 \\ \\
 & = \ds \sum_{n=0}^{t} v^n \sum_{\ell=0}^n (-1)^{\ell} \binom{n}{\ell}
 \binom{t -\beta (\ell +1)}{t}.
\end{array}
\eeq
The expression for the Sibuya state probabilities ${\mathbb P}(N_{Sibuya}(t)=n)$ and some related quantities were derived earlier \cite{PachonPolitoRicciuti2021}.
\begin{figure}[!t] 
\begin{center}
\includegraphics*[width=0.75\textwidth]{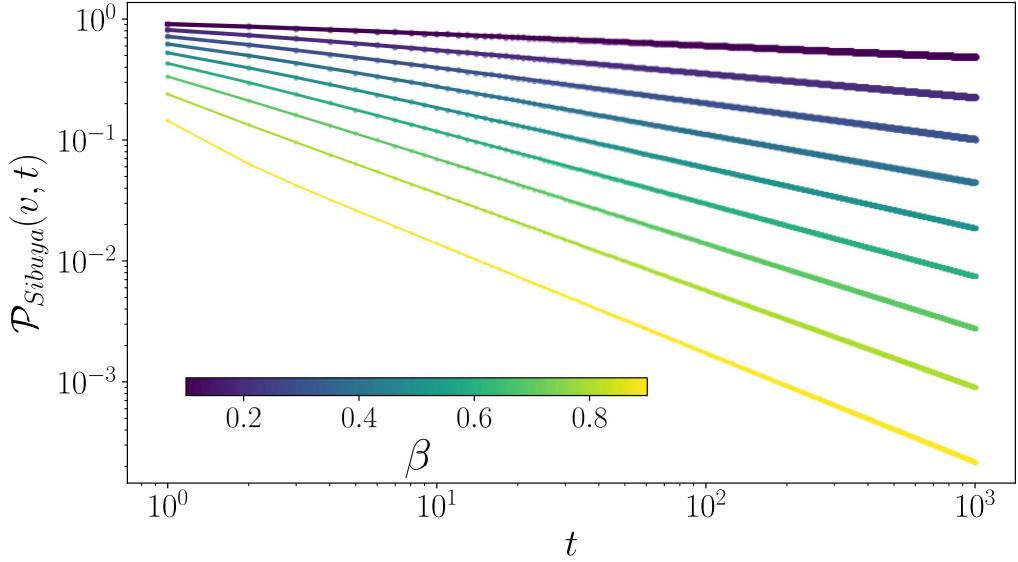}
\end{center}
\vspace{-5mm}
\caption{\label{Fig1}(Color online)
Sibuya state polynomial of Eq. (\ref{Sibuya-state_eval}) for $v=0.1$ and different values of $\beta \in (0,1)$ with extremely long waiting times for small $\beta$. The fat-tailed power-law decay emerges for large $t$ and can be seen in the slopes of the log-log representation.}
\end{figure}
In Figure \ref{Fig1} we plot the state polynomial for different values of $\beta$ where the non-markovianity with long-time memory of the Sibuya generator process is reflected by very long waiting times for small $\beta$ and shorter waiting times for larger $\beta$. One can identify the FT asymptotic power-law decay for large $t$ by the slopes, see also Eq. (\ref{asymptotic_state_poly}).
Then we have
\beq
\label{Lambda_sibuya}
\Lambda_{Sibuya}(a,b,t) =   \mathbb{E} a^{N_{Sibuya}(t)} b^{t- N_{Sibuya}(t)} = b^t{\cal P}_{Sibuya}\left(\frac{a}{b},t\right).
\eeq
The Sibuya ADTRW transition matrix then, with Eqs.
(\ref{evolution_eq}) and (\ref{Sibuya-state_eval}), yields
\beq
\label{Sibuya_trans}
{\mathbf P}_{Sib}(t) =\frac{(-1)^t}{t!} \sum_{n=0}^t({\mathbf W}^{-})^{t-n}({\mathbf W}^{+})^n
\sum_{\ell=0}^n (-1)^{\ell} \binom{n}{\ell} 
 \frac{\Gamma(\beta[\ell+1])}{\Gamma(\beta[\ell+1] -t)} ,\hspace{1cm} t \in \mathbb{N}_0
 \eeq
with the initial condition ${\mathbf P}_{Sibuya}(t)\big|_{t=0}= {\mathbf 1}$.
The terms for $\ell=0$ are all identical, namely
$\frac{\Gamma(t+1-\beta)} {\Gamma(1-\beta)\Gamma(t+1)}$ and dominating for $t \to \infty$. From Tauberian arguments it follows that this contribution is obtained from the dominating order for $u \to 1$ in the generating function,
namely the weakly singular term
${\bar \Phi}^{(n)}(u) \sim {\bar \Phi}^{(0)}(u)=(1-u)^{\beta-1}$. From this we see that all state probabilities have the same universal asymptotic scaling as the survival probability,
\beq
\label{survival_scaling}
{\mathbb P}(N_{Sibuya}(t)=n) \sim 
\lim_{t\to \infty} \frac{\Gamma(t+1-\beta)} {\Gamma(1-\beta)\Gamma(t+1)} =
\frac{t^{-\beta}}{\Gamma(1-\beta)} ,\hspace{0.5cm} \beta \in (0,1], \hspace{0.5cm} \forall n \in \mathbb{N}_0
\eeq
independent of $n$. This type of power-law scaling is universal for all fat-tailed waiting time distributions and it is equal to the asymptotic scaling of the Mittag-Leffler survival
probability, see relation (\ref{fat-tail_survival}). This fact can be generally
attributed to non-Markovianity and long-time memory \cite{GorenfloMainardi2008,Gorenflo2009} (consult also \cite{MichelitschPolitoRiascos2021}).
It follows that for the long-time asymptotic behavior of the state polynomial which then becomes an infinite power series in $v$, we have
\beq
\label{asymptotic_state_poly}
 {\cal P}_{Sibuya}(v,t) \sim \frac{1}{(1-v)} \frac{t^{-\beta}}{\Gamma(1-\beta)} , \hspace{1cm}
 \beta \in (0,1] , \hspace{0.5cm} t \to \infty.
 \eeq
For $\beta \to 0+$ extremely long-waiting times
occur between Sibuya successes (positive jumps), i.e.\
$\frac{t^{-\beta}}{\Gamma(1-\beta)} \to 1$. The limit $\beta \to 0$ corresponds to the above discussed 
frozen limit with very long waiting times between positive jumps and therefore strong domination of negative jumps.
For $\beta \to 1-$ the waiting times between Sibuya successes reaches their lower bound one, thus a strictly increasing markovian walk emerges where almost surely solely positive jumps occur, drawn from ${\mathbf W}^{+}$. In this Markovian limit the asymptotic relation  $\lim_{\beta \to 1-}\frac{t^{-\beta}}{\Gamma(1-\beta)} = \delta(t) =0 $ for $t$ large reflects the loss of memory.
\subsection{Simple Sibuya ADTRW}
\label{simple_Sibuya_ADRRW}
We now consider the simple Sibuya ADTRW, where almost surely only directed jumps of size one occur.
The transition matrix (\ref{simple_walk_transmat}) seen by the moving particle then yields

\beq
\label{trans_mat_sibuya_biased_walk}
\begin{array}{clr}
Q^{Sibuya}_{0,r}(t) & =  \ds  \frac{(-1)^t}{t!}\sum_{n=0}^{t} \delta_{0,r-2n}
\sum_{\ell=0}^n \binom{n}{\ell}(-1)^{\ell}  
 \frac{\Gamma(\beta[\ell+1])}{\Gamma(\beta[\ell+1] -t)},  & r \in \mathbb{Z} \\ \\
 &  =  \ds \frac{(-1)^t}{t!}
  \delta_{r,2\ceil{\frac{r-1}{2}}} \sum_{\ell=0}^{\frac{r}{2}}
  \binom{\frac{r}{2}}{\ell}(-1)^{\ell}  
 \frac{\Gamma(\beta[\ell+1])}{\Gamma(\beta[\ell+1] -t)}   &
  \end{array}
\eeq
supported on $r \in \{0,2,\ldots,2t-2,2t\}$. The Sibuya transition matrix then is related to Eq. (\ref{trans_mat_sibuya_biased_walk}) by
$P^{Sibuya}_{0,r}(t)=Q^{Sibuya}_{0,r+t}(t)$, see relations (\ref{Q_t_trans_mat}), (\ref{simple_walk_transmat}).
The Sibuya transition matrix solves with relation (\ref{simple_walk_renewal}) the renewal equations
\beq
\label{renewal_Sibuya}
P^{Sibuya}_{i,j}(t) = (-1)^t \binom{\beta-1}{t} \delta_{i,j+t}
+ \sum_{r=1}^t(-1)^{r-1}\binom{\beta}{r} P^{Sibuya}_{i,j+r-2}(t-r) ,\hspace{0.5cm} t \in \mathbb{N}
\eeq
with $P^{Sibuya}_{ij}(t)\big|_{t=0}=\delta_{ij}$
and those for $Q^{Sibuya}_{i,j}(t)= P^{Sibuya}_{i,j-t}(t)$ write by shifting the Kronecker symbols $\delta_{k,l} \to \delta_{k,l-t}$ in Eq. (\ref{renewal_Sibuya}).
The return probability to the departure site then is non-zero only for even $t=2s$
($s \in \{0,1,2, \ldots\}$) to give
\beq
\label{return}
\begin{array}{clr}
\ds P^{Sibuya}_{0,0}(t) & = \ds \delta_{t,\ceil{\frac{t-1}{2}}}\mathbb{P}[N_{Sibuya}(t)=\frac{t}{2}]  & \\ \\
& =  \ds \left\{ \begin{array}{l} \ds \frac{1}{t!} \sum_{\ell=0}^{\frac{t}{2}}
\binom{\frac{t}{2}}{\ell}(-1)^{\ell} 
\frac{\Gamma(\beta[\ell+1])}{\Gamma(\beta[\ell+1] -t)}
, \hspace{0.5cm} t\in \{0,2,4,\ldots\}  \\ \\
 0  , \hspace{1cm} t\in \{1,3,5,\ldots \}. \end{array}\right. &
 \end{array}
\eeq
\begin{figure}[!t] 
\begin{center}
\includegraphics*[width=0.75\textwidth]{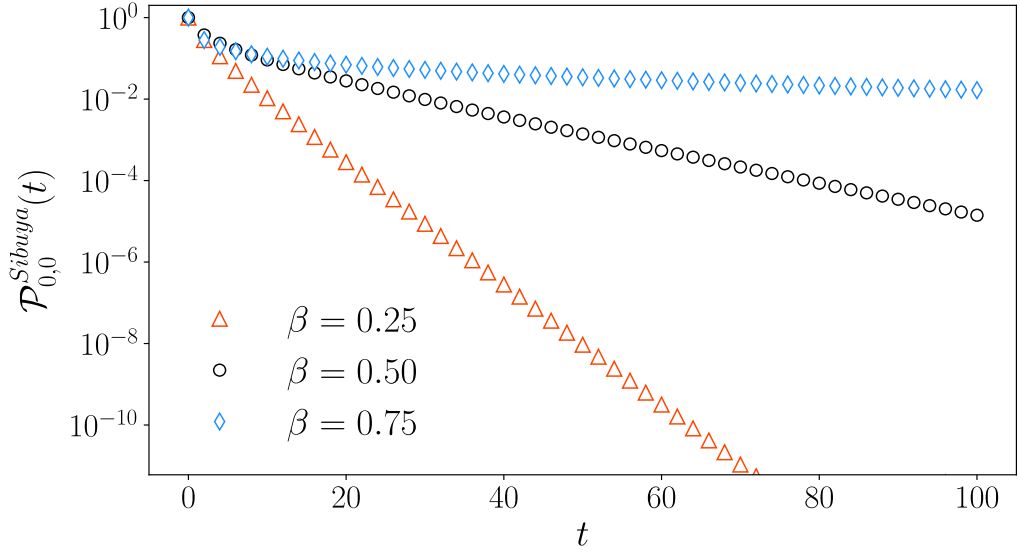}
\end{center}
\vspace{-5mm}
\caption{\label{Fig2}(Color online)
Return probability to the departure site of Eq.\ (\ref{return}) for different values of $\beta \in (0,1)$.}
\end{figure}
We plot in Figure \ref{Fig2} the Sibuya return probabilities to the departure site (\ref{return}) for different values of $\beta$. The smaller $\beta$ the more negative jumps $-1$ occur which reduces the return probability to the departure site thus the transient nature with a strong left-sided bias of the walk becomes more pronounced which can be clearly seen in Figure \ref{Fig3} depicting the expected position of the walker.
Figure \ref{Fig2} shows the tendency of the frozen limit $\beta \to 0+$ when negative jumps strongly dominate 
where the return probability to the departure site drops for $t>0$ `immediately' to zero.
\\[1mm]
For the issue of recurrence/transience indeed the generating function of the return probabilities
is especially important. In fact, this generating function is the diagonal element (see Eq.\,(\ref{complex_int}))
\beq
[{\bar \Lambda}_{Sibuya}({\hat T}_1,{\hat T}_{-1},u)]_{0,0} = \sum_{s=0}^{\infty} u^{2s} P^{Sibuya}_{0,0}(2s) = \sum_{s=0}^{\infty}
\frac{u^{2s}}{(2s)!} \sum_{\ell=0}^{s}
\binom{s}{\ell}(-1)^{\ell}
\frac{\Gamma(\beta[\ell+1])}{\Gamma(\beta[\ell+1] -2s)}
\eeq
and its canonical representation (see Eq.\,(\ref{gen_trans_feature}))
\beq
\label{canonic_diag}
\begin{array}{clr}
\ds [{\bar \Lambda}_{Sibuya}({\hat T}_1,{\hat T}_{-1},u)]_{0,0} & = \ds
\frac{1}{2\pi}
\int_{-\pi}^{\pi} {\bar \Lambda}_{Sibuya}(e^{-i\varphi},e^{i\varphi},u) {\rm d}\varphi & \\ \\  & = \ds \frac{1}{\pi} \Re \int_0^{\pi} \frac{e^{i\varphi}
 (1- u e^{i\varphi})^{\beta-1}}{2i\sin{\varphi} +
 e^{-i\varphi}(1-u e^{i\varphi})^{\beta}}{\rm d}\varphi. &
 \end{array}
\eeq
The EST on the departure site in an infinitely long walk yields (see also relation (\ref{tau_sojourn}))
\beq
\label{means_sojourn}
\mathbb{E}[ \tau_{00}]_{Sibuya} = [{\bar \Lambda}_{Sibuya}({\hat T}_1,{\hat T}_{-1},u)]_{0,0}\big|_{u=1} = \sum_{s=0}^{\infty}
\frac{1}{(2s)!} \sum_{\ell=0}^{s}
\binom{s}{\ell}(-1)^{\ell}
\frac{\Gamma(\beta[\ell+1])}{\Gamma(\beta[\ell+1] -2s)}.
\eeq
To explore this quantity it is convenient
to consider integral (\ref{canonic_diag}) at $u=1$, namely
\beq
\label{canonic_sojourn}
\ds \mathbb{E}[ \tau_{00}]_{Sibuya} =  \frac{1}{\pi} \Re \int_0^{\pi} \frac{e^{i\varphi}
 (1-e^{i\varphi})^{\beta-1}}{2i\sin{\varphi} +
 e^{-i\varphi}(1-e^{i\varphi})^{\beta}}{\rm d}\varphi.
\eeq
Taking into account that for $\varphi \to 0+$ we have
\beq
\label{small_varphi}
\Re {\bar \Lambda}_{Sibuya}(e^{-i\varphi},e^{i\varphi},1) \sim \Re \,\,
\frac{(-i)^{\beta-1}\varphi^{\beta-1}}{2i\varphi +
(-i)^{\beta}\varphi^{\beta}} \sim 2 \varphi^{-\beta} \cos\left(\frac{\pi\beta}{2}\right)
\eeq
being at $\varphi=0$ weakly singular as $\varphi^{-\beta}$ ($\beta \in (0,1)$). Hence $ {\bar \Lambda}_{Sibuya}(e^{-i\varphi},e^{i\varphi},1)$ is integrable at $\varphi=0$, by accounting for $A_{\beta}=1$, in agreement with the general FT feature (\ref{real_part}).
As a consequence the EST on the departure site (relations (\ref{means_sojourn}), (\ref{canonic_sojourn})) is finite. Therefore, the simple Sibuya ADTRW indeed is transient in agreement with our general proof for simple ADTRWs with FT waiting-time densities (case (a) in Section \ref{recurrence_transience}).
\\[1mm]
Let us now explore the bias and consider Eq. (\ref{simple_walk_bias}) by the generating function of the expected number of Sibuya arrivals, namely
\beq
\label{sibuya_arrivals}
\begin{array}{clr}
\ds {\bar N}_{Sibuya}(u) & =\ds  \frac{\partial }{\partial v} {\bar {\cal  P}}_{Sibuya}(v,u)\big|_{v=1} = \ds
\frac{1}{(1-u)}\frac{{\bar \psi}_{Sibuya}(u)}{(1-{\bar \psi}_{Sibuya}(u))}& \\ \\
 & = \ds (1-u)^{-\beta-1} -(1-u)^{-1} , & |u|<1
\end{array}
\eeq
yielding the exact non-negative expression 
\beq
\label{Sibuya_mean_arrival}
\begin{array}{clr}
\mathbb{E} N_{Sibuya}(t) & =  \ds  (-1)^t \binom{-(\beta+1)}{t} -1 = \frac{(\beta+1)\ldots (\beta +t)}{t!} -1 , & \hspace{1cm} t \in \mathbb{N} \\ \\
& = \ds \frac{\Gamma(\beta+t+1)}{\Gamma(\beta+1)\Gamma(t+1)} -1 =\binom{\beta + t}{t} -1 , & t \in \mathbb{N}_0
\end{array}
\eeq
holding for $\beta \in (0,1]$, see
also Ref. \cite{PachonPolitoRicciuti2021}.
The second line includes $t=0$ and reflects the initial condition $N_{Sibuya}(0)=0$.
We see in (\ref{Sibuya_mean_arrival}) that
$\mathbb{E} N_{Sibuya}(t)_{\beta \to 0+} \to 0 $ which means that for small $\beta$ the interarrival times between the Sibuya successes becomes infinitely long corresponding to the `frozen' limit. On the other hand we have the Markovian limit $\mathbb{E} N_{Sibuya}(t)_{\beta \to 1-} \to t$ where the trivial strictly increasing trivial walk emerges, see Section \ref{generator_process}.
With Eqs. (\ref{Sibuya_mean_arrival}), (\ref{simple_walk_bias}) we obtain for the expected position of the walker the exact expression
\beq
\label{Sibuya_simple_walk}
\mathbb{E}[ Y_t]_{Sibuya}  =  -t+   2 \mathbb{E} N_{Sibuya}(t) = 2\frac{\Gamma(\beta+t+1)}{\Gamma(\beta+1)\Gamma(t+1)} -2 -t ,\hspace{0.5cm} t \in \mathbb{N}_0, \hspace{1cm} \beta \in (0,1]
\eeq
where this result also holds in the Markovian limit $\beta=1$ and yields the upper
bound $\mathbb{E}[ Y_t]_{\beta=1}=t$. We also see that $\mathbb{E}[ Y_0 ]_{Sibuya} =0$ reflects the initial condition.
\begin{figure}[!t] 
\begin{center}
\includegraphics*[width=0.75\textwidth]{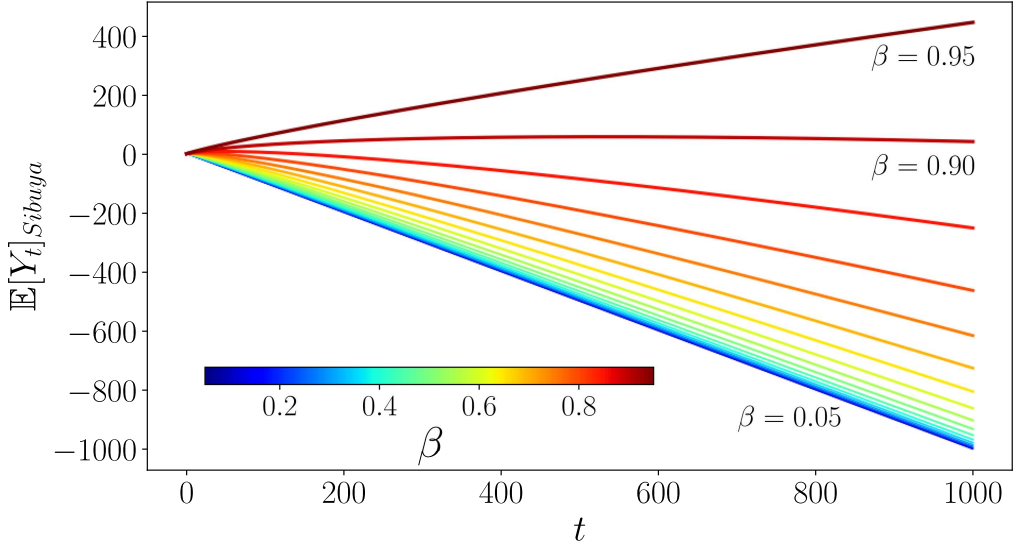}
\end{center}
\vspace{-5mm}
\caption{\label{Fig3}(Color online)
Expected position of the simple Sibuya ADTRW walker of Eq. (\ref{Sibuya_simple_walk}) for  different values of $\beta$. }
\end{figure}
In Figure \ref{Fig3} it is depicted the expected position $\mathbb{E}[ Y_t]_{Sibuya}$ for different values of $\beta \in (0,1)$. The smaller $\beta$ the longer the waiting times between positive jumps $+1$, the closer the 
lower bound $-t$ is approached. For $\beta \to 0+$ the `frozen' limit emerges
and the expected position approaches the lower bound $\lim_{\beta \to 0+}\mathbb{E}[ Y_t]_{Sibuya} = -t$ of a strictly decreasing walk with almost surely
solely jumps $-1$, see Section  \ref{generator_process}. On the other hand, for $\beta \to 1-$ where positive jumps $+1$ dominate one can see in the plot the larger $\beta$ the more the expected position approaches the upper bound $t$ reflecting $\lim_{\beta \to 1-} \mathbb{E}[ Y_t]_{Sibuya} = t$.
\\[1mm]
Consider now the asymptotic behavior with $\frac{\Gamma(\gamma+t)}{\Gamma(t)} \sim t ^{\gamma}$ for $t$ large. We have then the power-law
\beq
\label{mean-arrivals_tlarge}
\mathbb{E} [N]_{Sibuya}(t) \sim \frac{t^{\beta}}{\Gamma(\beta+1)}, \hspace{1cm} t \to \infty , \hspace{1cm} \beta \in (0,1]
\eeq
which is well known in anomalous diffusion (see e.g.\  \cite{GorenfloMainardi2008,Gorenflo2009,Laskin2003,MetzlerKlafter2001,MetzlerKlafter2004,MichelitschPolitoRiascos2021}) reflecting the non-Markovianity and
long memory feature of the Sibuya trial process. This relation includes the Markovian and frozen limits, respectively.
\\[1mm]
To see more closely the strong tendency of occurrence of negative jumps for the smaller $\beta$, as visible in Figure \ref{Fig3}, we consider the generating function of the expected position (\ref{Sibuya_simple_walk}) which takes, with Eq. (\ref{sibuya_arrivals}), the form
\beq\label{mean_gen}
  \sum_{t=0}^{\infty} \mathbb{E}[ Y_t]_{Sibuya} u^t =  2(1-u)^{-\beta-1} -2(1-u)^{-1}
  -u(1-u)^{-2} < 0 ,\hspace{1cm} u \to 1-.
  \eeq
We see in this relation the asymptotic behavior emerging in (\ref{Sibuya_simple_walk}) for $t \to \infty$, namely
\beq
\label{as_large_t}
\begin{array}{clr}
\mathbb{E}[ Y_t]_{Sibuya} &\sim 2 \frac{t^{\beta}}{\Gamma(\beta+1)} - t  \sim -t \to -\infty, \hspace{1cm} \beta \in (0,1) & \\ \\
 \mathbb{E}[ Y_t]_{Sibuya}& = t ,\hspace{2cm} \beta=1. & 
 \end{array}
\eeq
In the non-Markovian range $\beta \in (0,1)$ the long-time limit is governed by the term $-t$ due to domination of occurrence of negative jumps. Solely in the markovian limit $\beta=1$ the walk is strictly increasing with jumps $+1$ almost surely where the linear increase of the expected position
is maintained for all $t \in \mathbb{N}_0$.
\section{Asymmetric continuous-time random walk}
\label{time_changed_walk}
In the present section we introduce time-changed versions of the ADTRW.
We subordinate the ADTRW defined in Eq. (\ref{stepn}) to an independent renewal process, i.e.\ a continuous-time counting process $M(t) \in \mathbb{N}_0$ ($t \in \mathbb{R}^+$) with IID interarrival times such as Poisson, fractional Poisson and others.
We call the so defined walk `asymmetric continuous time random walk' (ACTRW). It turns out that the ACTRW is different from the classical Montroll--Weiss CTRW apart of some special cases also discussed in this section. ACTRWs are the class of random walks defined by
\beq
\label{time-changed}
{\cal Y}(t) = Y_{M(t)} = \sum_{j=1}^{M(t)} X_j ,\hspace{0.5cm} Y_0=0,
\hspace{0.5cm} X_j \in \mathbb{Z} \setminus \{0\}, \hspace{0.5cm} t \in \mathbb{R}^{+}
\eeq
where $Y_{m \in \mathbb{N}_0}$ is the ADTRW defined in Eq. (\ref{stepn}) with transition matrix (\ref{evolution_eq}).
In the ACTRW the trials of the generator process selecting the direction of the jumps $X_j$ take place
at the instants of arrival times of the point process $M(t)$.
The variable counting the arrivals in the composed process $N[M(t)] \in \mathbb{N}_0$ ($t\in \mathbb{R}^{+}$) 
indicates the number
of successes (number of positive jumps) in $M(t)$ trials 
and $M(t)-N[M(t)]$ the number of fails (number of negative jumps)
occurring within the continuous time interval $[0,t]$. 
The instants of successes are the 
continuous arrival times of the composed process $N[M(t)]$ which therefore is also a point process.
Compositions of counting processes (mainly of point processes)
where extensively studied in the literature \cite{OrsingherPolito2011,OrsingherPolito2012}.
\\[1mm]
Denoting with
$\mathbb{P}(M(t)=m)$ ($m\in \mathbb{N}_0$) the state probabilities (probabilities for $m$ arrivals within $[0,t]$) in the continuous-time process $M(t)$,
the state probabilities of the composed counting process $N[M(t)]$, i.e. the probabilities for $n$ arrivals (successes in the picture of trial process) within $[0,t]$) are given by
\beq
\label{state_prob_composition}
\mathbb{P}(N[M(t)] = n) =
\sum_{m=0}^{\infty} \mathbb{P}(M(t)=m)\mathbb{P}(N(m)=n) , \hspace{0.5cm} n \in \mathbb{N}_0, \hspace{0.5cm} t \in \mathbb{R}^{+}
\eeq
where we maintained for our convenience the vanishing
terms $m<n$ for which the state probabilities $\mathbb{P}(N(m)=n)=0$.
We call the composed continuous-time counting process $N[M(t)]$ the `time-changed generator process' of the ACTRW
since it contains information on the asymmetry of the walk: $N[M(t)]$ counts the number of positive jumps and $M(t)-N[M(t)]$ the number of negative jumps within $[0,t]$. 
The ACTRW transition matrix is then the time-changed version of Eq. (\ref{write_4_trans_mat}) and writes
\beq
\label{new_walk}
\begin{array}{clr}
\ds {\mathbf \Pi}({\mathbf W}^{+},{\mathbf W}^{-},t) & =  \ds \sum_{m=0}^{\infty} \mathbb{P}(M(t)=m)\sum_{n=0}^m\mathbb{P}(N(m)=n)  [{\mathbf W}^{+}]^n [{\mathbf W}^{-}]^{m-n} , & \hspace{0.5cm}t \in \mathbb{R}^{+}    \\ \\
\Pi_{i,j}(t)  & = \ds  \sum_{m=0}^{\infty} \mathbb{P}(M(t)=m) P_{i,j}(m)   
  \end{array}
\eeq
where ${\mathbf P}(m)=(P_{i,j}(m)) $ is the ADTRW transition matrix (\ref{write_4_trans_mat}). From the initial condition
$P_{i,j}(0)=\delta_{i,j}$ it follows $\Pi_{i,j}(t)\big|_{t=0}=\delta_{i,j}$, as a consequence of 
$\mathbb{P}(M(t)=m)\big|_{t=0}= \delta_{m,0}$.
The elements $\Pi_{i,j}(t)$ of the ACTRW transition matrix represent the probability that the walker is present on node $j$ at time $t$ with the indicated initial condition. Clearly, the ACTRW transition matrix (\ref{new_walk}) preserves the circulant property $\Pi_{i,j}(t)=\Pi_{0,j-i}(t)$. Since $M(t)$ is a continuous-time counting process, the random variable (\ref{time-changed}) describes a continuous-time random walk which is intrinsically asymmetric and - as we will see a little later - not of Montroll--Weiss type.
Then, it is useful to consider the scalar version of transition matrix (\ref{new_walk}),
\beq
\label{time-change}
\begin{array}{clr}
\ds  \Pi(a,b,t) & = 
 \ds \sum_{m=0}^{\infty}\mathbb{P}(M(t)=m) \sum_{n=0}^m\mathbb{P}(N(m)=n) a^n b^{m-n}, & \hspace{0.5cm} |a|,|b| \leq 1, \hspace{0.5cm} t \in \mathbb{R}^{+} \\ \\
  & = \ds \sum_{m=0}^{\infty}\mathbb{P}(M(t)=m) \Lambda(a,b,m),  &
 \end{array}
 \eeq
where $\Pi(1,1,t)=1$ reflects
the normalization of the state probabilities $\sum_{m=0}^{\infty}\mathbb{P}(M(t)=m)=1 $
with ${\cal P}(1,m)=\Lambda(1,1,m)=1$ (see normalization Eq. (\ref{normalization_state})) and
$\Pi(v,1,t)$ is the time-changed state polynomial.
This equation for $\Pi(a,b,t)$ shows the main difference to the Montroll--Weiss CTRW: For $b\neq 1$
the function $\Pi(a,b,t)$ and therefore the ACTRW transition matrix (\ref{new_walk}) are not represented by a series of the state probabilities (\ref{state_prob_composition}) of the composed process.
Therefore, the ACTRW generally is not in the Montroll--Weiss sense a random walk subordinated to the composed counting process $N[M(t)]$ (apart of some special cases such as the limits (\ref{Montroll_Weiss_case_i}), (\ref{Montroll_Weiss_case_ii}) and the example considered at the end of this section).
The general class of ACTRWs which can be reduced to Montroll--Weiss CTRWs have transition matrices 
of the form 
$\Lambda({\mathbf W},{\mathbf 1},t)$ where ${\mathbf W}$ is a single-jump transition matrix. The scalar version of this class
is obtained for $b=1$ in (\ref{time-change}) leading to (see
Eq. (\ref{state_prob_composition}))
\beq
\label{b_equal-to-one}
 \Pi(v,1,t)= \sum_{n=0}^{\infty} v^n \mathbb{P}(N[M(t)] = n) = \mathbb{E}_{N[M(t)]} v^{N[M(t)]}
\eeq
and is the time-changed state polynomial where the connection with the classical Montroll--Weiss CTRW can also be seen by means of its time-Laplace transform (\ref{time_changed_statepoly}).
\\[1mm]
For our further analysis it is convenient to consider the function (\ref{time-change}) in the Laplace domain.
The time-Laplace transform of a causal function $f(t)$ supported on $t \in \mathbb{R}^{+}$ is defined as
\beq
\label{Laplace_trafo_def}
{\tilde f}(s) =({\cal L} f) (s) = \int_0^{\infty}e^{-st}f(t){\rm d}t
\eeq
with a suitably chosen Laplace variable $s$. The time-Laplace transforms of the state probabilities $P(M(t)=m)$ reduce to
\beq
\label{laplace_state_innerprocess}
\int_0^{\infty}P(M(t)=m) e^{-st}{\rm d}t = \frac{1-{\tilde \eta}(s)}{s}({\tilde \eta}(s))^m ,\hspace{1cm} m \in \mathbb{N}_0
\eeq
where ${\tilde \eta}(s)$ denotes the Laplace transform of the interarrival time density $\eta(t)$ of the point process $M(t)$. The Laplace transform of function (\ref{time-change}) can hence be written as
\beq
\label{laplace_Pi}
\begin{array}{clr}
\ds {\tilde \Pi}(a,b,s) & = \ds \frac{1-{\tilde \eta}(s)}{s}\sum_{m=0}^{\infty}
({\tilde \eta}(s))^m b^m {\cal P}\left(\frac{a}{b},m\right) & |a|,|b| \leq 1 \\ \\
&= \ds
\frac{1-{\tilde \eta}(s)}{s} {\bar {\cal P}}\left(\frac{a}{b},\, b{\tilde \eta}(s)\right) = \frac{1-{\tilde \eta}(s)}{s} {\bar \Lambda}(a,b,{\tilde \eta}(s)) & \\ \\ & = \ds
\frac{1-{\tilde \eta}(s)}{s} \frac{b(1-{\bar \psi}[b {\tilde \eta}(s)])}{(1- b {\tilde \eta}(s) )(b-a{\bar \psi}[b {\tilde \eta}(s)])} &
\end{array}
\eeq
where necessarily ${\tilde \Pi}(1,1,s) = 1/s$ as a consequence of $\Pi(1,1,t)=1$ (reflecting the normalization condition of the state probabilities $\mathbb{P}(N[M(t)]=n)$ of the composed process).
In Eq. (\ref{laplace_Pi}) appears the generating function of the state polynomial ${\bar \Lambda}(a,b,{\tilde \eta}(s))= {\bar {\cal P}}(a/b,\, b{\tilde \eta}(s))$ (\ref{similar_gen}) with argument $u\to {\tilde \eta}(s)$ (fulfilling $|{\tilde \eta}(s)|\leq 1$), and ${\tilde \Pi}(v,1,s)$ is the Laplace transform of the time-changed state polynomial (\ref{b_equal-to-one}) having the simpler form
\beq
\label{time_changed_statepoly}
{\tilde \Pi}(v,1,s) =   \frac{1-{\bar \psi}[{\tilde \eta}(s)]}{s} \frac{1}{1-v{\bar \psi}[{\tilde \eta}(s)]} , \hspace{1cm} |v|\leq 1 .
\eeq
In this relation it appears the Laplace transform $\bar \psi[{\tilde \eta}(s)]$ of the waiting time density of the composed counting process $N[M(t)]$.
Therefore, Eq. (\ref{time_changed_statepoly}) has an interesting interpretation.
It is the time-Laplace transform of the generating function of the state probabilities ${\mathbb P}[N(M(t)) = n]$ ($n\in \mathbb{N}_0$) of the composed counting process. The interarrival time density of the composed process $N[M(t)]$ then reads
\beq
\label{assoc_iter_den}
\chi(t)= {\cal L}^{-1}\{{\bar \psi}[{\tilde \eta}(s)]\}(t) =\sum_{r=1}^{\infty}\psi(r) [\eta \star]^r(t) ,\hspace{0.5cm} t \in \mathbb{R}^{+}
\eeq
where we denote the inverse Laplace transform with
${\cal L}^{-1}\{\ldots\}(t)$. It follows from ${\bar \psi}[{\tilde \eta}(s)])\big|_{s=0}= {\bar \psi}(1)=1$ that $\chi(t)$ is indeed a density.
Then we have the Laplace transform of the state probabilities
(\ref{state_prob_composition}) as
\beq
\label{Laplace_trafo-comp}
\int_0^{\infty} e^{-st} \mathbb{P}(N[M(t)]=n)\,{\rm d}t = \frac{1-{\bar \psi}[{\tilde \eta}(s)]}{s} ({\bar \psi}[{\tilde \eta}(s)])^n,
\eeq
consistent with (\ref{time_changed_statepoly}).
In view of (\ref{laplace_Pi}) Laplace transform of the ACTRW transition matrix (\ref{new_walk}) can be written as
\beq
\label{laplace_trial_walk}
\ds  {\tilde \Pi}({\mathbf W}^{+},{\mathbf W}^{-},s) =
\frac{1-{\tilde \eta}(s)}{s} \frac{{\mathbf W}^{-}({\mathbf 1}-{\bar \psi}
[{\mathbf W}^{-}{\tilde \eta}(s)])}{({\mathbf 1}- {\mathbf W}^{-} {\tilde \eta}(s) )({\mathbf W}^{-}-{\mathbf W}^{+}{\bar \psi}[{\mathbf W}^{-} {\tilde \eta}(s)])}.
\eeq
Clearly, this expression
does not have a Montroll--Weiss structure.
However, in some special cases, for instance for the limits $\alpha_k \leq \epsilon \to 0+$ and $\alpha_k \to 1-$ $\forall k$, respectively strictly decreasing and increasing Montroll--Weiss CTRWs emerge, namely
\beq
\label{Montroll_Weiss_case_i}
\begin{array}{clr}
\ds  \Pi_{0+}(t) & =  \ds \sum_{m=0}^{\infty}\mathbb{P}(M(t)=m)({\mathbf W}^{-})^m ,& \hspace{0.5cm} \alpha_k \leq \epsilon   \to 0+ \\ \\
\ds {\tilde \Pi}_{0+}(s)  & =  \ds  \frac{(1-{\tilde \eta}(s))}{s} [{\mathbf 1}-{\mathbf W}^{-}{\tilde \eta}(s))]^{-1} &
\end{array}
\eeq
and
\beq
\label{Montroll_Weiss_case_ii}
\begin{array}{clr}
\ds   \Pi_{1-}(t) & = \ds \sum_{m=0}^{\infty}\mathbb{P}(M(t)=m)({\mathbf W}^{+})^m  , & \hspace{0.5cm} \alpha_t \to 1- \\ \\
\ds  {\tilde \Pi}_{1-}(s) & = \ds  \frac{(1-{\tilde \eta}(s))}{s} [{\mathbf 1}-{\mathbf W}^{+}{\tilde \eta}(s))]^{-1}. &
\end{array}
\eeq
\subsection{ACTRW evolution equations}
With the above considerations we can derive the time-evolution equations for $\Pi(a,b,t)$ and the
ACTRW transition matrix.
To this end we rearrange Eq. (\ref{laplace_Pi}) to
\beq
\label{rearrange_laplace}
 \frac{1-{\bar \psi}[b{\tilde \eta}(s)]}{s{\bar \psi}[b{\tilde \eta}(s)]}\left[
 s {\tilde \Pi}(a,b,s) -1  \right] = -\frac{(1-b){\tilde \eta}(s)}{1-b{\tilde \eta}(s)}\frac{(1-{\bar \psi}[b{\tilde \eta}(s)])}{s{\bar \psi}[b{\tilde \eta}(s) ]}+
\left(\frac{a}{b}-1\right){\tilde \Pi}(a,b,s).
\eeq
Introducing the auxiliary kernels
\beq
\label{auxiliary_kernel}
{\cal K}(b,t) = {\cal L}^{-1}\left\{ \frac{1-{\bar \psi}[b{\tilde \eta}(s)]}{s{\bar \psi}[b{\tilde \eta}(s)]} \right\}(t)
,\hspace{0.5cm} t \in \mathbb{R}^{+}, \hspace{0.5cm} |b| \leq 1
\eeq
and
\beq
\label{axiliary_kernel2}
{\cal R}(b,t) = {\cal L}^{-1} \left\{ \frac{{\tilde \eta}(s)}{1-b{\tilde \eta}(s)} \right\}(t)
 \eeq
where we observe that $\frac{(1-b){\tilde \eta}(s)}{1-b{\tilde \eta}(s)}\big|_{s=0}=1$ since ${\tilde \eta}(s)\big|_{s=0}=1$ thus $(1-b){\cal R}(b,t)$ is a normalized density.
${\cal K}(1,t)$ can be seen as the memory kernel of the composed process $N[M(t)]$.
Writing Eq. (\ref{rearrange_laplace}) in the time-domain yields the following Cauchy problem
\beq
\label{time_domain_evolution}
\begin{array}{clr}
\ds \int_0^t {\cal K}(b,t-\tau)\frac{\mathrm d}{\mathrm d\tau} \Pi(a,b,\tau)\,{\rm d}\tau & = \ds
(b-1)\int_0^t {\cal K}(b,t-\tau){\cal R}(b,\tau)\,{\rm d}\tau + \left(\frac{a}{b}-1\right)\Pi(a,b,t), & \hspace{0.5cm} t \in \mathbb{R}^+\\ \\
\ds  \Pi(a,b,t)\big|_{t=0}   = 1 & &
\end{array}
\eeq
and defines also the Cauchy problem for the ACTRW transition matrix (\ref{new_walk}) by replacing $a$ with ${\mathbf W}^{+}$, $b$ with ${\mathbf W}^{-}$ and considering the initial condition $\Pi_{ij}|_{t=0}=\delta_{ij}$. 
The left-hand side of the scalar equation (\ref{time_domain_evolution}) is a general fractional derivative and has profound connections to the general fractional calculus introduced by Kochubei \cite{Kochubei2011} and see also \cite{Giusti2020,Luchko2021,Diethelm2020,TMM_APR_2020,TMM-FP-APR-Gen-Mittag-Leffler2020}.
The influence of the asymmetry can be seen in the change of sign
of the second term on the right-hand side for $a>b$ and $a< b$ for real $a,b$, respectively. We also recover for $a=b=1$ that the right hand side 
is null thus $\Pi(1,1,t)=1$ is constant (the normalization of the state probabilities 
(\ref{state_prob_composition}) as a conserved quantity $\forall t $ as $a=b$ introduces a further symmetry --- one should recall the previously mentioned connection with Noether's theorem).
The difference to the Montroll--Weiss CTRW becomes obvious by the presence of the first term on the right-hand side for $b\neq 1$.
For $b=1$ we have $(1-b){\cal R}(b,t)\big|_{b=1}=0$ and
the auxiliary kernel (\ref{auxiliary_kernel}) reduces to the memory kernel of the composed counting process $N[M(t)]$.
Thus, Eq.\,(\ref{time_domain_evolution}) then reduces to the form of a generalized Kolmogorov--Feller equation
of Montroll--Weiss type 
for $\Pi(v,1,t)$, being solved by the Laplace-inverse of Eq.\,(\ref{time_changed_statepoly}).
\\[1mm]
Eq.\,(\ref{time_domain_evolution}) governs the time-evolution in a ACTRW and is the counterpart to the generalized Kolmogorov--Feller equation which occurs in a Montroll--Weiss CTRW. 
\\[3mm]
As a pertinent example let us consider the point process  $M_{\mu}(t)$ to be the time-fractional Poisson process and the independent trial process $N_B(m)$ to be the Bernoulli process. We will see that in contrast to the general case, in this example the ACTRW indeed boils down to a Montroll--Weiss type CTRW.
The Laplace transform of the waiting time density of the time-fractional Poisson process $M_{\mu}(t)$ has the form ${\tilde \eta}_{\mu}(s)=\frac{\xi_0}{\xi_0+s^{\mu}}$ ($\xi_0>0, \mu \in (0,1]$)  \cite{Laskin2003} and the Bernoulli waiting time generating function is ${\bar \psi}_B(z)=pz/(1-qz)$ ($p+q=1$), thus the waiting time density of the composed process $N_B[M_{\mu}(t)]$ has the Laplace transform 
\beq
\label{Laplace_fract-poisson_bernoulli}
{\bar \psi}_B[{\tilde \eta}_{\mu}(s)] = \frac{p\xi_0}{p\xi_0+s^{\mu}} , \hspace{1cm} \mu \in (0,1]
\eeq
i.e.\ the composition is also a continuous time fractional Poisson process with changed 
constant $\xi=p\xi_0$. If $\mu=1$ the standard Poisson process is recovered.
The auxiliary kernel (\ref{auxiliary_kernel}) yields
\beq
\label{auxiliary-FracBern}
{\cal K}_{\mu}(b,t) = {\cal L}^{-1}\left(\frac{s^{\mu-1}}{p\xi_0b} +\frac{1-b}{bp}s^{-1} \right)  = \frac{ t^{-\mu}}{bp\xi_0\Gamma(1-\mu)}+ \frac{1-b}{bp} , \hspace{1cm} t \geq 0
\eeq
and the second kernel (\ref{axiliary_kernel2})
is expressed by a Mittag--Leffler density
\beq
\label{we_get_for-the_second}
(1-b){\cal R}_{\mu}(b,t)  = 
{\cal L}^{-1}\left(\frac{\xi_0(1-b)}{\xi_0(1-b)+s^{\mu}} \right)  =
\xi_0(1-b)t^{\mu-1} E_{\mu,\mu}(-\xi_0(1-b)t^{\mu})=-\frac{d}{dt}E_{\mu}(-\xi_0(1-b)) ,\hspace{0.5cm} t\geq 0
\eeq
where $E_{\mu,\gamma}(z)$ and $E_{\mu}(z)$ denote the generalized and the standard Mittag--Leffler functions, respectively (see \cite{SamkoKilbasMarichev1993,OldhamSpanier1974} for definitions and properties). We introduce the Caputo-fractional derivative \cite{SamkoKilbasMarichev1993,OldhamSpanier1974}
\beq
\label{Caputo_fract-der}
\frac{d^{\mu}}{dt^{\mu}} y(t) = \int_0^t\frac{ (t-\tau)^{-\mu}}{\Gamma(1-\mu)}\frac{d}{d\tau}y(\tau){\rm d}\tau ,\hspace{1cm} \mu \in (0,1]
\eeq
recovering for $\mu\to 1-$ the standard first order derivative. The Cauchy problem (\ref{time_domain_evolution}) after some routine manipulations takes the form of a fractional differential equation
\beq
\label{fractional_Cauchy}
\begin{array}{clr}
\ds \frac{d^{\mu}}{dt^{\mu}}\Pi_{\mu,\lambda}(a,b,t)
&=\ds  -\lambda \, \Pi_{\mu,\lambda}(a,b,t) , & \ds  \hspace{1cm} \lambda = \xi_0[1- qb - p a)] , \hspace{0.5cm} \mu \in (0,1]\\ \\
\ds \Pi_{\mu,\lambda}(a,b,t)\bigg|_{t=0} & =\ds  1 &
\end{array}
\eeq
where keep also in mind that $|a|,|b| \leq 1$. Cauchy problem
(\ref{fractional_Cauchy}) then has the Mittag-Leffler solution
\beq
\label{solution_ACTRW_frac-ber}
\Pi_{\mu,\lambda}(a,b,t) ={\cal L}^{-1}\left(\frac{s^{\mu-1}}{s^{\mu}+\lambda}\right)   = E_{\mu}(-\lambda t^{\mu}) 
\eeq
which is also directly obtained from Eq. (\ref{laplace_Pi}).
The ACTRW with the time-changed generator process $N_B[M_{\mu}(t)]$ has therefore the Mittag--Leffler transition-matrix
\beq
\label{ML-Jumps}
{\mathbf P}_{\mu}(t)= \Pi_{\mu,\lambda}({\mathbf W}^{+},{\mathbf W}^{-},t) = E_{\mu}(-\xi_0t^{\mu}[{\mathbf 1}-p{\mathbf W}^{+}-q{\mathbf W}^{-}]).
\eeq
One can see by means of the Laplace transforms
that this transition matrix has the particularity that it is of Montroll--Weiss type where Bernoulli jumps with a well defined `Laplacian matrix' ${\mathbf 1}-p{\mathbf W}^{+}-q{\mathbf W}^{-}$ are subordinated to the independent time-fractional Poisson process $M_{\mu}(t)$.
%where for $\mu=1$ the Poisson case is recovered where the Mittag-Leffler functions take their exponential counterparts.
%
%
%
%
\section{Conclusions}
\label{Conclusions}
We have presented a new type of asymmetric discrete-time random walk, the ADTRW. In this walk the direction of the jumps is determined by the outcomes of a trial process (the `generator process') which is constructed as a discrete-time counting process. We considered the ADTRW on the integer line and analyzed recurrence/transience features. We demonstrated that fat-tailed waiting time distributions in the generator process generate transience and bias in a simple ADTRW whereas light-tailed waiting-time distributions allow both transient and recurrent behavior. In the recurrent case the simple ADTRW is unbiased in an asymptotic sense (i.e.\ in the limit $t\to \infty$).
We proved that among the simple ADTRWs solely the one with symmetric Bernoulli generator process is strictly unbiased in the sense that the expected position is null (i.e.\ on the departure site) at all times.
On the other hand for all transient cases
the simple ADTRW is biased and vice versa. The ADTRW model can be generalized to several directions. For instance modifications in the trial process for the determination of directions of the jumps define new types of ADTRWs with a large potential of new applications. Further possible generalizations include involvement of long-range jumps where different one-jump transition matrices are selected by counting processes, or the possibility of further considering the
interarrival times between two consecutive fails (negative jumps) not geometrically distributed.
\\[1mm]
We also considered prescribed admissible functions for the expected position $\mathbb{E}[ Y_t]_{simple}= f(t)$ in a simple ADTRW, see Eq. (\ref{leads_to}).
For future research interesting candidates are constituted by the class of discrete-time versions of non-negative Bernstein functions which are strictly positive $\mathbb{E}[ Y_t]_{simple} = f(t) >0$ for $t\in \mathbb{N}$. 
The special interest of this topic is also due to the possibility to construct simple ADTRWs which in the Ruin Game interpretation provide strategies where the expectation value of the assets never hits the ruin condition 
(zero assets).
For an analytical procedure 
to construct discrete approximations of Bernstein functions, see \cite{TMM-FP-APR-Gen-Mittag-Leffler2020} and consult also \cite{PachonPolitoRicciuti2021,MichelitschPolitoRiascos2021}.
\\[1mm]
We also introduced time-changed versions of the ADTRW leading to the ACTRW model. The ACTRW constitutes a new class of biased continuous-time random walks which are generally not of Montroll--Weiss type, apart of some special cases. 
In the present paper we could only introduce the main idea of the ACTRW model which merits further thorough analysis and exploration of pertinent cases. 
\\[1mm]
The new types of asymmetric random walks introduced in the present paper open a wide field of interdisciplinary applications in `complex systems' such as in finance, birth and death models, and biased anomalous transport and diffusion.

\begin{appendix}
\section{\small APPENDICES}
\subsection{Some pertinent limits of $\Lambda(a,b,t)$}
\label{Appendix_A}
We consider here some limiting cases of the state polynomial $\Lambda(a,b,t)$ defined in Eq. (\ref{similar}). 
\\[1mm]
First let $b\to 0$ for which we get for the related generating function
\beq\label{limit_b_0}
\begin{array}{clr}
 \bar{\Lambda}(a,0,u) & =  \ds  \lim_{b \to 0} {\bar {\cal  P}}\left(\frac{a}{b},ub\right)
=  \lim_{b \to 0} \frac{1}{1-\frac{a}{b}{\bar \psi}(bu)}
=\sum_{t=0}^{\infty} \mathbb{P}(N(t)=t) u^t a^t  , \hspace{0.5cm} \mathbb{P}(N(t)=t)=\alpha_1^t &  \\ \\  & = \ds  \frac{1}{1-a\alpha_1 u}. &
  \end{array}
\eeq
Thus 
\beq
\label{markovian_lim}
\Lambda(a,0,t) = \mathbb{P}(N(t)=t) a^t = (a\alpha_1)^t
\eeq
containing only the order $t$ of the state polynomial where we account for the probability of $t$ successes in $t$ trials
$\mathbb{P}(N(t)=t)=\alpha_1^t$
and $ \lim_{b \to 0}  {\bar \psi}(bu)/b =u\alpha_1$.
\\[1mm]
A further pertinent limit is obtained for $a=0$, namely
\beq
\label{limit_a_zero}
  {\bar \Lambda}(0,b,u) = {\bar {\cal  P}}\left(\frac{a}{b},ub\right)\Big|_{a=0} =
 \frac{1-{\bar \psi(bu)}}{1-{\bar bu}} = \sum_{t=0}^{\infty} (bu)^t
 \mathbb{P}(N(t)=0)  
\eeq
retrieving the (rescaled) survival probability
\beq
\label{frozen-lim}
\Lambda(0,b,t) =  b^t\mathbb{P}(N(t)=0).
\eeq
Plainly, these limiting relations are connected with the `Markovian' and `frozen' limits, respectively, see Section \ref{generator_process}. 
\subsection{Time shift operator representations}
\label{Appendix_B}
In order to derive some convenient operator representations of above deduced renewal and master equations
(\ref{renewal_eq_Lam})-(\ref{master_eq_memory}) be reminded that we deal with causal
(discrete-time) distributions supported on non-negative integers
\beq
\label{causal}
F(t) =\Theta(t) f(t) , \hspace{1cm} t \in \mathbb{Z}
\eeq
i.e.\ they are null for negative $t$ which we indicate by the discrete Heaviside function $\Theta(t)$ defined in Eq. (\ref{discrete-theta}).
Then, we introduce the time-backward shift operator
${\hat {\cal T}}_{-1}$ which is such that ${\hat {\cal T}}_{-1}f(t)=f(t-1)$, where $t \in \mathbb{Z}$ is the discrete time coordinate. We further use throughout the paper the following equivalence of generating functions and shift
operators (see \cite{MichelitschPolitoRiascos2021} for an outline of essential properties)
\beq
\begin{array}{clr}
\ds f(t) & = \ds \frac{1}{t!} \frac{d^t}{du^t}{\bar f}(u)|_{u=0} & \hspace{0.5cm} t \in \mathbb{N}_0 \\ \\
& = \ds\sum_{k=0}^{\infty} f(k) {\cal T}_{-k} \delta_{0,t},  & \hspace{0.5cm} {\cal T}_{-k} \delta_{0,t} = \delta_{0,t-k} =\delta_{k,t} \\ \\
  & = \ds {\bar f}({\hat {\cal T}}_{-1}) \delta_{0,t}. &
\end{array}
\eeq
We see that causality of a distribution $f(t)$ is generated by ${\bar f}({\hat {\cal T}}_{-1})\delta_{0,t}$, where in the generating function we replaced $u$ with ${\hat {\cal T}}_{-1}$ and only non-negative powers 
of the time backward shift operator ${\bar T}_{-1}$ are considered (${\cal T}_{-1}^k={\cal T}_{-k}$).
The interarrival time density has then the 
backward time shift operator representation
\beq
\label{time_shift_interarrival_density}
\psi(t) = {\bar \psi}({\cal T}_{-1})\delta_{0,t} = \sum_{k=1}^{\infty} \psi(k) {\cal T}_{-k}\delta_{0,t}
= \sum_{k=1}^{\infty} \psi(k) \delta_{0,t-k}.
\eeq
The state polynomial can be represented as
\beq
\label{solved_by}
{\cal P}(v,t)={\bar P}(v,{\hat {\cal T}}_{-1})\delta_{0,t} =\frac{1-{\bar \psi}({\hat {\cal T}}_{-1})} { [1- {\hat {\cal T}}_{-1}][1-v{\bar \psi}({\hat {\cal T}}_{-1})]} \delta_{0,t}
\eeq
and 
\beq
\label{Lambda_shift_representation}
\Lambda(a,b,t)=b^t{\cal P}\left(\ds \frac{a}{b},t\right)  = \frac{1-{\bar \psi}(b{\hat {\cal T}}_{-1})}{1-b{\hat {\cal T}}_{-1}}  \frac{1}{1-\frac{a}{b} {\bar \psi}(b{\hat {\cal T}}_{-1})} \delta_{0,t} , \hspace{1cm} b\neq 0 .
\eeq
Note that shift operators and shift operator functions commute among each other reflecting the commutative property of discrete convolutions.
Then by accounting for ${\bar f}(b {\hat {\cal T}}_{-1})\delta_{0,t} = b^t {\bar f}({\hat {\cal T}}_{-1})\delta_{t,0}= b^tf(t)$ we can rewrite this relation as
\beq
\label{represent_Lambda_2}
\begin{array}{clr}
\ds \Lambda(a,b,t) &  =  \ds
\frac{1-{\bar \psi}(b{\hat {\cal T}}_{-1})}{1-b{\hat {\cal T}}_{-1}} \delta_{0,t}
+  
\frac{a}{b} {\bar \psi}(b{\hat {\cal T}}_{-1}) \Lambda(a,b,t) & \\ \\ &  =  \ds b^t\Phi^{(0)}(t) + \frac{a}{b} {\bar \psi}(b{\hat {\cal T}}_{-1}) \Lambda(a,b,t) &
\end{array}
\eeq
and use that
\beq
\label{convolution_generfunction}
\begin{array}{clr}
\ds {\bar \psi}(b{\hat {\cal T}}_{-1}) \Lambda(a,b,t)
& = \ds  \sum_{r=1}^{\infty} b^r \psi(r) {\hat {\cal T}}_{-r} \Lambda(a,b,t)  & \\ \\
& =  \ds \sum_{r=1}^{\infty} b^r  \psi(r) \Lambda(a,b,t-r) = \sum_{r=1}^{t} b^r  \psi(r) \Lambda(a,b,t-r)
\end{array}
\eeq
where in the last line we used causality, i.e.  $\Lambda(a,b,t-r)=0$ when $r>t$. 
We hence arrive at the renewal equation (\ref{renewal_eq_Lam}), i.e.
\beq
\label{end-result_Lambda_renewal}
 \Lambda(a,b,t) =  b^t\Phi^{(0)}(t) + \sum_{r=1}^{t} a b^{r-1}  \psi(r) \Lambda(a,b,t-r).
\eeq
\subsection{Bernoulli trials}
\label{Appendix_C}
As the simplest and best known special case we briefly recall some well-known features of the memoryless Bernoulli walk. This walk matches in our ADTRW model as `simple Bernoulli ADTRW' where
$\alpha_t =p$ $\forall t \in \mathbb{N}$ does not depend on $t$, leading to geometric waiting time density (\ref{define_renewal}) $\psi_B(t) =pq^{t-1}$ ($p+q=1$, $t\in \mathbb{N}$), with generating function ${\bar \psi}_B(u)=\frac{pu}{1-qu}$ and survival probability $\Phi^{(0)}_B(t)=q^t$ ($t\in \mathbb{N}_0$). It is straight-forward to see that the Bernoulli state probabilities
are given by the {\it Binomial distribution}
\beq
\label{Binomial_distribution}
\begin{array}{clr}
\ds \Phi_B^{(n)}(t) & = \ds \frac{1}{t!}\frac{d^t}{du^t} \frac{p^nu^n}{(1-qu)^{n+1}}\bigg|_{u=0}= \frac{t!}{n!(t-n)!}p^nq^{t-n} , \hspace{1cm} n \leq t &\\ \\  & & \ds  n, t \in \mathbb{N}_0.  \\
 & = \ds  0 , \hspace{1cm}  n > t  &
 \end{array}
\eeq
The state polynomial (\ref{PZ}) yields then straight-forwardly
\beq
\label{state-Bernoulli}
{\cal P}_B(v,t) = (pv+q)^t 
\eeq
and
\beq
\label{Lambda_bern}
\Lambda_B(a,b,t) = \mathbb{E} \, a^{N(t)} b^{t-N(t)}=(pa+qp)^t.
\eeq
This relation contains information on the bias where the expected position of the walker $\mathbb{E}[ (Y_B)_t] $ at time $t$ in a simple walk (\ref{simple_walk_bias}) (i.e.\ with unit next neighbor jumps) leads to the well-known classical result
 \cite{Feller1971,RednerS}
\beq
\label{mean}
  \mathbb{E}[ (Y_B)_t]  = \left(\frac{\partial}{\partial a} -\frac{\partial}{\partial b}\right) \Lambda_B(a,b,t)\bigg|_{a=b=1} =  (p-q)t.
\eeq
A measure for the asymmetry of the walk is provided here by ${\cal C}_{B} = t^{-1}\mathbb{E}[ (Y_B)_t] =p-q$ where for $p=q=\frac{1}{2}$ this simple walk is strictly unbiased.
The EST (\ref{tau_sojourn}) then writes (with transition matrix ${\mathbf P}_B(t)= [p{\mathbf W}^{+}+q{\mathbf W}^{-}]^t$)
\beq
\label{bernoulli-meansojounr}
\mathbb{E}[ \tau_{rs} ]_B = \sum_{t=0}^{\infty}[\left(p{\mathbf W}^{+}+q{\mathbf W}^{-}\right)^t]_{rs} =
\left(\left[{\mathbf 1}-p{\mathbf W}^{+}-q{\mathbf W}^{-}\right]^{-1}\right)_{rs} .
\eeq
For more details consult 
\cite{SpitzerF1976,Feller1971,PolyaG1921,RednerS,LawlerLimick2012}, and the references therein.

\subsection{Some features of light-tailed waiting time densities}
\label{Appendix_complex} 
We introduce the auxiliary generating function

\beq
\label{gdis}
{\bar g}(z)= \sum_{t=1}^{\infty}\psi(t)z^{t-1} ,\hspace{1cm} {\bar \psi}(z)= z{\bar g}(z)
\eeq
with ${\bar g}(0)=\psi(1)=\alpha_1$ and ${\bar g}(1)={\bar \psi}(1) =1$ and where ${\bar \psi}(z)$ denotes the generating function (\ref{genfu}) of a discrete-time waiting time density supported on $\mathbb{N}$ (where $0<\psi(1)=\alpha_1\leq 1$ and $\alpha_1=1$ only in the trivial case when ${\bar \psi}_{trivial}(z)=z$ with constant ${\bar g}_{trivial}(z)=1$ and $\psi_{trivial}(t)=\delta_{t1}$ with unit waiting times). We consider here only the situation in which $\psi(t)$ is LT.
Then we introduce
\beq
\label{prop1}
g_1 = \frac{d}{dz}{\bar g}(z)\big|_{z=1} = 
\sum_{t=2}^{\infty}(t-1)\psi(t) =A_1-1 \geq 0  ,\hspace{1cm} z\in [-1,1]
\eeq
where $A_1= \frac{d}{dz}{\bar \psi}(z)\big|_{z=1}=\sum_{t=1}^{\infty}\psi(t)t \geq 1$ denotes the expected waiting time. 
For some related properties of generating functions, we refer to the book of Harris \cite{Harris_book}.
The aim of this appendix is to prove for LT waiting-time densities (with $ \psi(t)=g(t-1)$, i.e.\ ${\bar \psi}(z)=z{\bar g}(z)$) the existence of the canonical representation
\beq
\label{canonic_function-rep}
{\bar g}(z)-z = (z-1)(z-r) e^{h(z)} ,\hspace{1cm} |z|\leq 1
\eeq
at least on the unit disc where $h(z)$ is analytic. The zero $r$ is real and non-negative.
The convex property of ${\bar g}(z)$ with $g_1=\frac{d}{dz}{\bar g}(z)\big|_{z=1}=A_1-1 \geq 0$ allows us to determine the properties of the zero $r=r(A_1) \in \mathbb{R}^{+}$: 
\beq
\label{properties_zeros}
\begin{array}{clr}
\ds  r(A_1)  >   1 , & \hspace{0.5cm} A_1 <2 & \\ \\
\ds r(A_1) = 1, & \hspace{0.5cm} A_1=2 & \\ \\
\ds 0 < r(A_1) <1 ,  & \hspace{0.5cm} A_1> 2. &
\end{array}
\eeq
Note that $r$ is outside of the unit disc for $A_1<2$, inside for $A_1>2$, and $r=1$ for $A_1=2$. 
Since $r(A_1)$ is a continuous function of $A_1$ we can infer that
\begin{itemize}
\item \noindent $\lim_{A_1\to 1+0} r(A_1) = + \infty$ (limit of trivial walk);
\item \noindent  $\lim_{A_1 \to \infty} r(A_1) = 0 $ (FT limit).
\end{itemize}
Further limiting properties can be seen from relations (\ref{properties_zeros}) together with the continuity of
$r(A_1)$. 
\\[1mm]
This behavior can be explicitly seen in the Bernoulli trial process: see Eq.
(\ref{poles}) where $r_B=\frac{p}{q}$ and $A_1=\frac{1}{p}$ thus ${\bar g}_{B}(z)-z= (z-1)(z-r_B)e^{h_B(z)}$ with $h_B(z)=\log(q)-\log(1-qz)$.
\\[1mm]
To prove the existence of canonical representation (\ref{canonic_function-rep}) let us now consider the cases $1\leq A_1<2$ and $A_1>2$ separately.
\\[3mm]
{\bf Case} $1\leq A_1<2$ ($g_1=\frac{d}{dz}{\bar g}(z)\big|_{z=1} \in [0,1)$):
\\
One observes then the following properties:
\\
\noindent (i) \hspace{3cm} $ \ds \frac{d^{n}}{dz^n}{\bar g}(z) \geq 0 ,\hspace{1cm} z\in [-1,1], \hspace{1cm} n\in \mathbb{N}_0$
\\ \\i.e.\ ${\bar g}(z)$ is absolutely monotonic (AM) for $ z\in [-1,1]$. Hence ${\bar g}(-z)$ is in that interval completely monotonic (CM) and strictly positive. Therefore we have
\\ \\ 
\noindent (ii) \hspace{3cm}  $\ds  0 \leq \frac{d}{dz}{\bar g}(z) \leq \frac{d}{dz}{\bar g}(z)\big|_{z=1}= g_1 = A_1-1 <1 , \hspace{1cm} z \in [-1,1]$ 
\\ \\ \\
\noindent (iii) \hspace{3cm} $\ds  0< {\bar g}(-1) \, < \,  {\bar g}(z) \, < {\bar g}(1) \, = \, 1 $ , \hspace{0.3cm} ${\bar g}(z)$ is AM for $z\in [-1,1]$.
\\ \\ 
In (ii) we have $g_1=0$ only for $A_1=1$, otherwise $0 < g_1$.
\\
\noindent Remark: There is a connection of ${\bar \psi}(z)$ with Bernstein functions, see 
\cite{TMM-FP-APR-Gen-Mittag-Leffler2020} (Formula (39)).
\\[3mm]
\noindent Remark to (i)-(iii): 
\\ 
The derivatives $\ds \frac{d^{n}}{dz^n}{\bar g}(z) = 
\sum_{t=n+1}^{\infty} \frac{(t-1)!}{(t-1-n)!}\psi(t)z^{t-1-n}  \geq 0$ are AM functions for
$z \in [0,1]$ and for all $n\in \mathbb{N}_0$. As a consequence ${\bar g}(-z)$ 
is completely monotonic (CM) with $(-1)^n\frac{d^{n}}{dz^n}{\bar g}(-z) \geq 0$ for $z\in [0,1]$.
Then, it follows that
${\bar g}(-z) = \sum_{t=1} \psi(t)(-z)^{t-1}$ is strictly positive on the interval $z \in [0,1]$ with the minimal value ${\bar g}(-1)$. In that interval
${\bar g}(-1) < {\bar g}(-z) < {\bar g}(0)=\alpha_1 \leq 1$. From this inequality, reflecting the AM feature of ${\bar g}(z)$, it follows the property (iii), i.e.\ ${\bar g}(z)$ is strictly positive and AM
on the real interval $z \in [-1,1]$.
From (i)--(iii) we see that 
\beq
\label{we_have}
\begin{array}{clr}
{\bar g}(z) &  >z ,\hspace{1.5cm} -1 \leq z <1 & \\ \\
               {\bar g}(z)    &  = z ,\hspace{1.5cm} z_1=1 &
               \end{array}
                   \eeq
i.e.\ $z_1=1$ is the only zero of the function (\ref{canonic_function-rep}) in the real interval $[-1,1]$ when $1\leq A_1<2$, therefore we have in that range $r(A_1)>1$.
\\[3mm]
Recall that in this part our goal is to prove the canonical representation (\ref{canonic_function-rep}) first for $g_1<1$ ($1\leq A_1<2$).
To this end we make use of 
\beq
\label{wemake_use}
\left|\frac{d}{dz}{\bar g}(|z|e^{i\varphi})\right| \leq \left|\frac{d}{dz}{\bar g}(|z|)\right| \leq g_1 <1 ,\hspace{1cm} |z| \leq 1 
\eeq
and therefore 
\beq
\label{has_the_rep}
\frac{d}{dz}{\bar g}(z) = a(x,y) e^{i\alpha(x,y)}
\eeq
where $a(x,y) = |\frac{d}{dz}{\bar g}(z)|$ denotes the absolute value and
$\alpha(x,y)$ the argument of $\frac{d}{dz}{\bar g}(z)$.
If there are complex zeros of the function (\ref{canonic_function-rep}) for $|z|<1$ then they solve $|y|= |\Im\{ {\bar g}(x+iy)\}|$ with $y=\Im(z) \neq 0$. Now consider the imaginary part 

\beq
\label{imaginary}
\begin{array}{clr}
\ds \Im\{ {\bar g}(x+iy) \} & = \ds \Im\left\{ {\bar g}(x)+ \int_0^y \frac{d}{dx}{\bar g}(x+i\tau)i{\rm d}\tau \right\} = \Im\left\{\int_0^y \frac{d}{dx}{\bar g}(x+i\tau)i{\rm d}\tau \right\} ,\hspace{0.5cm} |y|\leq |z|< 1 & \\ \\
 & = \ds \int_0^y a(x,\tau)  e^{i\alpha(x,\tau)}{\rm d}\tau, & 
 \end{array}
\eeq
having absolute value
\beq
\label{absolute_value}
\begin{array}{clr} 
|\Im\{ {\bar g}(x+iy)|  =   | \Im\int_0^y a(x,\tau)  e^{i\alpha(x,\tau)}{\rm d}\tau | \leq  | \int_0^y a(x,\tau)  e^{i\alpha(x,\tau)}{\rm d}\tau |
\leq \int_0^y a(x,\tau) {\rm d}\tau < & g_1 |y| < |y| , &
\end{array}
\eeq
i.e. the left-hand side is always smaller than $|y|$ as $a(x,y) < g_1=A_1-1<1$ where this inequality holds for any $|z|\leq 1$.
Therefore $y=\Im\{ {\bar g}(x+iy)\}$ has for $x\in [-1,1]$ only the trivial solution
$y=0$, which concludes the proof that
(\ref{canonic_function-rep}) for $A_1<2$ has solely the zero $z_1=1$ and no zeros for $|z|<1$.
We hence can establish the following theorem:
\\ 
{\bf Theorem I}\\ {\it Let ${\bar g}(x)$ be non-negative and absolutely monotonic (AM) and ${\bar g}(x) \neq  x$ on the real interval $x\in [0,\rho_0)$ with
$|\frac{d}{dz}{\bar g}(z)| \leq g_1(\rho_0) < 1$ for $|z| \leq \rho_0$. Then, the complex function ${\bar g}(z)-z \neq 0$ has no zeros within the disc $|z| < \rho_0$.}
\\ 
We will use this theorem to complete the 
proof of the canonical representation (\ref{canonic_function-rep}) for $A_1>2$.
\\[3mm]
{\bf Case} $A_1>2$ ($g_1= \frac{d}{dz}{\bar g}(z)\big|_{z=1} \in (1,\infty)$):
\\ 
From the convexity of ${\bar g}(x)$ for $x \in (0,1)$ it follows that there is a real zero $r=r(g_1)$ which is for $g_1>1$ in the interval $r(g_1) \in(0,1)$. We have then for $x \in [0,1)$,
\beq
\label{has_zero_g}
\begin{array}{clr}
{\bar g}(x)-x & > 0 ,\hspace{1cm} 0 \leq x < r  & \\[3mm]
  {\bar g}(r)-r & = 0  ,\hspace{1cm} r \in (0,1)  &      \\[3mm]
{\bar g}(x)-x & < 0 ,\hspace{1cm} r < x < 1. &
\end{array}
\eeq
The inequality of the first line can be seen from
${\bar g}(0)=\alpha_1 >0$ and using the continuity of ${\bar g}(x)$.
On the other hand we have $|{\bar g}(z)| \leq {\bar g}(|z|)$ ($z=x+iy$) as a consequence of the AM property and together with the third line of (\ref{has_zero_g}) the inequality
\beq
\label{first_int}
|{\bar g}(z)| \leq {\bar g}(|z|) < |z| ,\hspace{1cm}  |z| \in (r,1)
\eeq
i.e. there are no zeros within the ring $|z| \in (r,1)$
\beq
\label{no-zero_interval}
{\bar g}(z) - z \neq 0 , \hspace{1cm} |z| \in (r,1).
\eeq
Now we finally consider $|z|<r$. From the convexity of
${\bar g}(x)$ it follows that 
\beq
\label{derivative_in_r}
\left|\frac{d}{dz}{\bar g}(z)\right| <  \frac{d}{dz}{\bar g}(z)\Big|_{z=r} = g_1(r) < 1 ,\hspace{1cm} |z|<r
\eeq
i.e., the slope $\frac{d}{dz}{\bar g}(z)\big|_{z=r}<1$ must be smaller than one
in the cutting point ${\bar g}(x)=x$, that is at $x=r$.
Since $\frac{d}{dz}{\bar g}(z)\big|_{z=r}<1$ apply Theorem I with $\rho_0=r$
to see that
\beq
\label{no_complex_zero}
{\bar g}(z)-z \neq 0 , \hspace{1cm} |z|<r .
\eeq
The existence of the inequalities (\ref{no-zero_interval}), (\ref{no_complex_zero}) concludes the proof 
that for $g_1>1$ there are no
further zeros for $|z|\leq 1$ except the real zeros $z=r \in (0,1)$ and $z=1$.
\\[1mm]
In conclusion the canonical form (\ref{canonic_function-rep}) with the property (\ref{properties_zeros}) holds true at least for $|z|\leq 1$ in the admissible range $A_1 \in [1,\infty)$.
\\[3mm]
{\it Poisson waiting time density} 
\\ [1mm]
As a prototypical example for an LT waiting-time density we consider Poisson waiting times with
$\psi_P(\lambda,t)= e^{-\lambda} \frac{\lambda^{t-1}}{(t-1)!}$ ($t \in \mathbb{N}$) and 
${\bar \psi}_P(\lambda,z)=z{\bar g}_P(\lambda,z)$ with ${\bar g}_P(\lambda,z)= e^{\lambda(z-1)}$ where
$g_1 =\lambda = A_1-1$ and $\lambda \geq 0$ ($\lambda=0$ corresponding to the trivial walk). Clearly, ${\bar g}_P(\lambda,z)$ is AM on $[-1,1]$. Consider
inequality (\ref{absolute_value}) for $0\leq \lambda < 1$. It yields
\beq
\label{expon-zeros}
 |e^{\lambda(x-1)} \Im\{ e^{i\lambda y}\}| =  e^{\lambda(x-1)} |\sin{\lambda y}| < |y|.
\eeq
Therefore
\beq
\label{equality_for_zeros}
e^{\lambda(z-1)} -z =0 , \hspace{1cm} 0 \leq  \lambda <1
\eeq
has the only zero $z_1=1$ for $|z|\leq 1$. We also see that
\beq
{\bar g}_P(\lambda,z)-z= e^{\lambda(z-1)}-z = (z-1)\left[\lambda-1 +\sum_{\ell=2}^{\infty}\frac{\lambda^{\ell}}{\ell!} (z-1)^{\ell-1}\right]
\eeq
has zero $z_1=1$ with multiplicity one for $\lambda \neq 1$ and multiplicity
$2$ for $\lambda=1$ (i.e.\ $A_1=2$) which is the recurrent limit. We prove now that the canonical form 
\beq
\label{canonic_function-rep2}
e^{\lambda(z-1)}-z = (z-1)(z-r_{\lambda}) e^{h_{\lambda}(z)} , \hspace{1cm} z \in \mathbb{C}, \hspace{1cm} {\bar \psi}(z)= z e^{\lambda(z-1)}
\eeq
holds here for all $\lambda \geq 0$ ($A_1=\lambda+1$).
The zero $r_{\lambda} \in \mathbb{R}^{+}$ has the properties (\ref{properties_zeros}) where $g_1=\lambda$ ($A_1=\lambda+1$).
Let us explore whether
\beq
\label{z_intersection}
z=e^{\lambda(z-1)} , \hspace{1cm} \lambda > 0 
\eeq
has complex zeros with non-vanishing imaginary parts.
Applying the logarithm to both sides we have $\log(z)=\lambda(z-1)$ and
with $z=\rho e^{i\varphi}$, considering the real and imaginary parts, it yields the two equations
\beq
\label{real_and_imag}
\begin{array}{clr}
\log\rho & = \lambda(\rho\cos(\varphi)-1) & \\ \\
\varphi & = \lambda \rho \sin(\varphi) & 
\end{array}
\eeq
each defining a parametrization of lines $\rho_{1,2}(\varphi)$
where both sides of Eq. (\ref{z_intersection}) have the same argument $\varphi$. %
The intersections of Eq. (\ref{z_intersection}) are the points where $\rho_{1}(\varphi)=
\rho_2(\varphi)$. The second equation gives
\beq
\label{rho_2}
\rho_2(\varphi,\lambda)= \frac{\varphi}{\lambda \sin(\varphi)} ,\hspace{1cm} \varphi \neq 0
\eeq
Plugging Eq. (\ref{rho_2}) into the first equation of (\ref{real_and_imag}) we obtain
\beq
\label{rho_1}
\rho_1(\varphi,\lambda)=  e^{\varphi\cot(\varphi)-\lambda}.
\eeq
We will see that $\rho_2(\varphi,\lambda) -\rho_1(\varphi,\lambda) > 0$ by rewriting the inequality
\beq
\label{no_intersecton}
 \frac{\varphi}{\sin(\varphi)} > \lambda e^{-\lambda}  e^{\varphi\cot(\varphi)}
\eeq
holding for all $\lambda \geq 0$ and $\varphi \neq 0$ (which we show below) to conclude the proof of existence of canonical form (\ref{canonic_function-rep2}) with property (\ref{properties_zeros}). 
\\[3mm]
We verify inequality (\ref{no_intersecton}) as follows: $\frac{\varphi}{\sin(\varphi)}\big|_{\varphi \to 0}=1 > \lambda e^{-\lambda}$ $\forall \lambda \geq 0$. Thus this inequality is true in the limit $\varphi \to 0$.
Then consider the slope of the left hand side
\beq
\label{slope_left}
\frac{d}{d\varphi}\left(\frac{\varphi}{\sin(\varphi)}\right)= \frac{\sin(\varphi)-\varphi\cos(\varphi)}{\sin^2(\varphi)} ,\hspace{1cm} \varphi \in [-\pi,\pi]
\eeq 
This is an odd function with the same sign as $\varphi$, i.e.\ negative for $\varphi <0$, null for $\varphi=0$ and positive for $\varphi>0$, thus
$\frac{\varphi}{\sin(\varphi)} \geq 1$ and the value $1$ at $\varphi=0$ is a minimum
for $\varphi \in [-\pi,\pi]$.
Then consider the slope of the right-hand side of inequality (\ref{no_intersecton}),
\beq
\label{slope_right}
 \lambda e^{-\lambda} \frac{d}{d\varphi} \,  e^{\varphi \frac{\cos\varphi}{\sin\varphi}} =
   \lambda e^{-\lambda} e^{\varphi \frac{\cos\varphi}{\sin\varphi}}
  \frac{(\sin(2\varphi)-2\varphi)}{2\sin^2\varphi} 
 \eeq
which is an odd function with the opposite sign than $\varphi$, i.e.\ 
for $\varphi \neq 0$, i.e. the value at $\varphi=0$ is a maximum. Thus $\frac{\varphi}{\sin{\varphi}} \geq 1 >  \lambda e^{-\lambda}\geq  \lambda e^{-\lambda} e^{\varphi \frac{\cos\varphi}{\sin\varphi}}$ and hence inequality (\ref{no_intersecton}) holds true for $\varphi \in [-\pi,\pi]$ ($\varphi \neq 0$), i.e.\ Eq. (\ref{z_intersection}) has no intersections with non-vanishing imaginary parts. 
\\[3mm]
{\it Non-negativity of EST (\ref{sojournm1})}
\\
This formula holds for $A_1<2$ which we rewrite as

\beq
\label{sojournm-rewrite}
\mathbb{E}[\tau_{0,-1}]_{simple} =  \frac{1+g_1}{1-g_1}  - \frac{1}{\alpha_1} , \hspace{0.5cm} 0\leq g_1 <1
\eeq
where from Eq. (\ref{prop1}) follows that $g_1 \geq \sum_{t=2}^{\infty}\psi(t) = 1-\alpha_1$. Thus $\alpha_1\geq 1-g_1$ and so $\mathbb{E}[\tau_{0,-1}]_{simple} \geq 0 $ for all non-trivial cases $1-\alpha_1 \leq  g_1<1$ (and for the trivial walk ${\bar g}(z)=1$ is constant with $g_1=0$ and $\alpha_1=1$ thus $\mathbb{E}[\tau_{0,-1}]_{simple}=0$) concluding the proof.

\end{appendix}

\end{document}